\title{Pseudo-random graphs
}
\author{
Michael Krivelevich\thanks{
Department of Mathematics, Raymond and Beverly Sackler Faculty of Exact
Sciences, Tel Aviv University, Tel Aviv 69978, Israel. E-mail:
krivelev@post.tau.ac.il. Research supported in part by a USA-Israel
BSF Grant, by a grant from the Israel Science Foundation
 and by a Bergmann Memorial Grant.}
\and
Benny Sudakov \thanks{Department of Mathematics,
Princeton University, Princeton, NJ 08544, USA. 
Email address: bsudakov@math.princeton.edu.
Research supported in part by NSF grant DMS-0106589.
Part of this research was done while visiting
Microsoft Research.}                        
}
\documentclass[11pt]{article}
\usepackage{latexsym}
\usepackage{amsfonts}
\usepackage{amsmath}
\newtheorem{theo}{Theorem}[section]
\newtheorem{prop}[theo]{Proposition}

\newtheorem{coro}[theo]{Corollary}
\newtheorem{conj}[theo]{Conjecture}

\newcommand{\Epsilon}{{\cal E}}
\def\om{\omega}

\date{}
\begin{document}
\maketitle

\section{Introduction}

Random graphs have proven to be one of the most important and fruitful
concepts in modern Combinatorics and Theoretical Computer Science.
Besides being a fascinating study subject for their own sake, they
serve as essential instruments in proving an enormous number of
combinatorial statements, making their role quite hard to overestimate.
Their tremendous success serves as a natural motivation for the
following very general and deep informal questions: what are the
essential properties of random graphs? How can one tell when a given
graph behaves like a random graph? How to create deterministically
graphs that look random-like? This leads us to a concept of {\em
pseudo-random graphs}.

Speaking very informally, a pseudo-random graph $G=(V,E)$ is a graph
that behaves like a truly random graph $G(|V|,p)$ of the same edge
density $p=|E|\left/{{|V|}\choose 2}\right.$. Although the last sentence
gives some initial idea about this concept, it is not very informative,
as first of all it does not say in which aspect the pseudo-random graph
behavior is similar to that of the corresponding random graph, and
secondly it does not supply any quantitative measure of this
similarity. There are quite a few possible graph parameters that can
potentially serve for comparing pseudo-random and random graphs (and in
fact quite a few of them are equivalent in certain, very natural sense,
as we will see later), but probably the most important characteristics
of a truly random graph is its {\em edge distribution}. We can thus make
a significant step forward and say that a pseudo-random graph is a
graph with edge distribution resembling the one of a truly random graph
with the same edge density. Still, the quantitative measure of this
resemblance remains to be introduced.

Although first examples and applications of pseudo-random graphs 
appeared very long
time ago, it was Andrew Thomason who launched systematic research on
this subject with his two papers \cite{Tho87a}, \cite{Tho87b} in the
mid-eighties. Thomason introduced the notion of jumbled graphs, enabling
to measure in quantitative terms the similarity between the edge
distributions of pseudo-random and truly random graphs. He also supplied
several examples of pseudo-random graphs and discussed many of their
properties. Thomason's papers undoubtedly defined directions of future
research for many years.

Another cornerstone contribution belongs to Chung, Graham and Wilson
\cite{ChuGraWil89} who in 1989 showed that many properties of different
nature are in certain sense equivalent to the notion of
pseudo-randomness, defined using the edge distribution. This fundamental
result opened many new horizons by showing additional facets of
pseudo-randomness.

Last years brought many new and striking results on pseudo-randomness
by various researchers. There are two clear trends in recent research on
pseudo-random graphs. The first is to apply very diverse methods from
different fields (algebraic, linear algebraic, combinatorial,
probabilistic etc.) to construct and study pseudo-random graphs. The
second and equally encouraging is to find applications, in many cases
quite surprising, of pseudo-random graphs to problems in Graph Theory,
Computer Science and other disciplines. This mutually enriching
interplay has greatly contributed to significant progress in research on
pseudo-randomness achieved lately. 

The aim of this survey is to provide a systematic treatment of the
concept of pseudo-random graphs, probably the first since the two
seminal contributions of Thomason 
\cite{Tho87a}, \cite{Tho87b}. Research in pseudo-random graphs has
developed tremendously since then, making it impossible to provide full
coverage of this subject in a single paper. We are thus forced to omit
quite a few directions, approaches, theorem proofs from our discussion.
Nevertheless we will attempt to provide the reader with a rather
detailed and illustrative account of the current state of research in
pseudo-random graphs.

Although, as we will discuss later, there
are several possible formal approaches
to pseudo-randomness, we will mostly emphasize the approach based on
graph eigenvalues. We find this approach, combining linear algebraic
and combinatorial tools in a very elegant way, probably the most
appealing, convenient and yet quite powerful.

This survey is structured as follows. In the next section we will
discuss various formal definitions of the notion of pseudo-randomness,
from the so called jumbled graphs of Thomason to the $(n,d,\lambda)$-graphs defined by Alon,
where pseudo-randomness is connected to the eigenvalue gap. We then
describe several known constructions of pseudo-random graphs, serving both as
illustrative examples for the notion of pseudo-randomness, and also as
test cases for many of the theorems to be presented afterwards. The
strength of every abstract concept is best tested by properties it
enables to derive. Pseudo-random graphs are certainly not an exception
here, so in Section 4 we discuss various properties of pseudo-random
graphs. Section 5, the final section of the paper, is devoted to
concluding remarks.

\section{Definitions of pseudo-random graphs}
Pseudo-random graphs are much more of a general concept describing some
graph theoretic phenomenon than of a rigid well defined notion -- the
fact reflected already in the plural form of the title of this
section! Here we describe various formal approaches to the concept of
pseudo-randomness. We start with stating known facts on the edge
distribution of random graphs, that will serve later as a benchmark for
all other definitions. Then we discuss the notion of jumbled graphs
introduced by Thomason in the mid-eighties. Then we pass on to the
discussion of graph properties, equivalent in a weak (qualitative) sense
to the pseudo-random edge distribution, as revealed by Chung, Graham and
Wilson in \cite{ChuGraWil89}. Our next item in this section is the
definition of pseudo-randomness based on graph eigenvalues -- the
approach most frequently used in this survey. Finally, we discuss the
related notion of strongly regular graphs, their eigenvalues and their
relation to pseudo-randomness.

\subsection{Random graphs}
As we have already indicated in the Introduction, pseudo-random graphs are
modeled after truly random graphs, and therefore mastering the edge
distribution in random graphs can provide the most useful insight on
what can be expected from pseudo-random graphs. The aim of this
subsection is to state all necessary definitions and results on random
graphs. We certainly do not intend to be comprehensive here, instead
referring the reader to two monographs on random graphs \cite{Bol01},
\cite{JanLucRuc00}, devoted entirely to the subject and presenting a
very detailed picture of the current research in this area.

A {\em random graph} $G(n,p)$ is a probability space of all labeled
graphs on $n$ vertices $\{1,\ldots,n\}$, where for each pair 
$1\le i<j\le n$, 
$(i,j)$ is an edge of $G(n,p)$ with probability $p=p(n)$,
independently of any other edges. Equivalently, the probability of a
graph $G=(V,E)$ with $V=\{1,\ldots,n\}$ in $G(n,p)$ is
$Pr[G]=p^{|E(G)|}(1-p)^{{n\choose 2}-|E(G)|}$. We will occasionally
mention also the probability space $G_{n,d}$, this is the probability
space of all $d$-regular graphs on $n$ vertices endowed with the uniform
measure, see the  survey of Wormald \cite{Wor99} for more background.
We also say that a graph property ${\cal A}$ holds {\em almost surely}, or a.s. for
brevity, in $G(n,p)~(G_{n,d})$ if the probability that $G(n,p)~(G_{n,d})$ 
has ${\cal A}$ tends to one as the number of vertices $n$ tends to 
infinity.

From our point of view  the most important parameter of random graph
$G(n,p)$ is its edge distribution. This characteristics can be easily
handled due to the fact that $G(n,p)$ is a product probability
space with independent appearances of different edges. Below we cite
known results on the edge distribution in $G(n,p)$.

\begin{theo}\label{dist1}
Let $p=p(n)\le 0.99$. Then almost surely $G\in G(n,p)$ is such that if
$U$ is any set of $u$ vertices, then
$$
\left|e(U)-p{u\choose 2}\right| =
O\left(u^{3/2}p^{1/2}\log^{1/2}(2n/u)\right)\ .
$$
\end{theo}

\begin{theo}\label{dist2}
Let $p=p(n)\le 0.99$. Then almost surely $G\in G(n,p)$ is such that 
if $U,W$ are disjoint sets of vertices satisfying $u=|U|\le w=|W|$, then
$$
\left|e(U,W)-puw\right|=O\left(u^{1/2}wp^{1/2}\log^{1/2}(2n/w)\right)\ .
$$
\end{theo}

The proof of the above two statements is rather straightforward. Notice
that both quantities $e(U)$ and $e(U,W)$ are binomially distributed
random variables with parameters ${u\choose 2}$ and $p$, and $uw$ and
$p$, respectively. Applying standard Chernoff-type estimates on the
tails of the binomial distribution (see, e.g., Appendix A of
\cite{AloSpe00}) and then the union bound, one gets the desired
inequalities.

It is very instructive to notice that we get less and less control over
the edge distribution as the set size becomes smaller.
 For example, in the probability
space $G(n,1/2)$ every subset is expected to contain half of its
potential edges. While this is what happens almost surely for large
enough sets due to Theorem \ref{dist1}, there will be almost surely
sets of size about $2\log_2n$ containing all possible edges (i.e.
cliques), and there will be almost surely sets of about the same size, 
containing no edges at all (i.e. independent sets). 

For future comparison we formulate the above two theorems in the
following unified form:

\begin{coro}\label{dist3}
Let $p=p(n)\le 0.99$. Then almost surely in $G(n,p)$ for every two
(not necessarily) disjoint subsets of vertices $U,W\subset V$ of
cardinalities $|U|=u, |W|=w$, the number $e(U,W)$ of edges of $G$
with one endpoint in $U$ and the other one in $W$ satisfies:
\begin{eqnarray}\label{dist}
|e(U,W)-puw|= O(\sqrt{uwnp})\ .
\end{eqnarray}
\end{coro}
(A notational agreement here and later in the paper: if an edge $e$
belongs to the intersection $U\cap W$, then $e$ is counted twice in
$e(U,W)$.)

Similar bounds for edge distribution hold also in the space $G_{n,d}$
of $d$-regular graphs, although they are significantly harder to derive
there. 

Inequality (\ref{dist}) provides us with a quantitative benchmark,
according to which we will later measure the uniformity of edge 
distribution in pseudo-random graphs on $n$ vertices with edge density
$p=|E(G)|\left/{n\choose 2}\right.$

It is interesting to draw comparisons between research in random graphs
and in pseudo-random graphs. In general, many properties of random
graphs are much easier to study than the corresponding properties of
pseudo-random graphs, mainly due to the fact that along with the almost
uniform edge distribution described in Corollary \ref{dist3}, random
graphs possess as well many other nice features, first and foremost of
them being that they are in fact very simply defined product
probability spaces. Certain graph properties can be easily shown to hold
almost surely in $G(n,p)$ while they are not necessarily valid in
pseudo-random graphs of the same edge density. We will see quite a few
such examples in the next section. A general line of research
appears to be not to use pseudo-random methods to get new results for
random graphs, but rather to try to adapt techniques developed
for random graphs to the case of pseudo-random graphs, or alternatively
to develop original techniques and methods.

\subsection{Thomason's jumbled graphs}
In two fundamental papers \cite{Tho87a}, \cite{Tho87b} published in 1987
Andrew Thomason introduced the first formal quantitative definition of
pseudo-random graphs. It appears quite safe to attribute the launch of
the systematic study of pseudo-randomness to Thomason's papers. 

Thomason used the term "jumbled" graphs in his papers. A graph
$G=(V,E)$ is said to be $(p,\alpha)$-{\em jumbled} if $p,\alpha$ are
real numbers satisfying $0<p<1\le \alpha$ if every subset of vertices
$U\subset V$ satisfies:
\begin{eqnarray}\label{jumbled}
\left|e(U)- p{{|U|}\choose 2}\right|\le \alpha |U|\ .
\end{eqnarray}
The parameter $p$ can be thought of as the density of $G$, while
$\alpha$ controls the deviation from the ideal distribution. According
to Thomason, the word "jumbled" is intended to convey the fact that the
edges are evenly spread throughout the graph. 

The motivation for the above definition can be clearly traced to the
attempt to compare the edge distribution in a graph $G$ to that of a
truly random graph $G(n,p)$. Applying it indeed to $G(n,p)$ and
recalling (\ref{dist}) we conclude that the random graph $G(n,p)$ is
almost surely $O(\sqrt{np})$-jumbled.

Thomason's definition has several trivial yet very nice features.
Observe for example that if $G$ is $(p,\alpha)$-jumbled then the
complement $\bar{G}$ is $(1-p,\alpha)$-jumbled. Also, the definition is
hereditary -- if $G$ is $(p,\alpha)$-jumbled, then so is every induced
subgraph $H$ of $G$.

Note that being $(p,\Theta(np))$-jumbled for a graph $G$ on $n$
 vertices and ${n\choose 2}p$ edges does not say too much about the edge
distribution of $G$ as the number of edges in linear sized sets can
deviate by a percentage from their expected value.  However
as we shall see very soon if $G$ is known to be $(p,o(np))$-jumbled,
quite a lot can be said about its properties. Of course, the smaller is
the value of $\alpha$, the more uniform or jumbled is the edge
distribution of $G$. A natural question is then how small can be the
parameter $\alpha=\alpha(n,p)$ for a graph $G=(V,E)$ on $|V|=n$ vertices
with edge density $p=|E|\left/{n\choose 2}\right.$? Erd\H os and
Spencer proved in \cite{ErdSpe72} that $\alpha$ satisfies
$\alpha=\Omega(\sqrt{n})$ for a constant $p$; their method can be
extended to show $\alpha=\Omega(\sqrt{np})$ for all values of $p=p(n)$.
We thus may think about $(p,O(\sqrt{np}))$-jumbled graphs on $n$
vertices as in a sense best possible pseudo-random graphs.

Although the fact that $G$ is $(p,\alpha)$-jumbled carries in it a lot
of diverse information on the graph, it says almost nothing (directly at
least) about small subgraphs, i.e. those spanned by subsets $U$ of size
$|U|=o(\alpha/p)$. Therefore in principle a $(p,\alpha)$-jumbled graph
can have subsets of size $|U|=O(\alpha/p)$ spanning by a constant factor
less or more edges then predicted by the uniform distribution. In many
cases however quite a meaningful local information (such as the presence
of subgraphs of fixed size) can still be salvaged from global 
considerations as we will see later.

Condition (\ref{jumbled}) has obviously a global nature as it applies to
{\em all} subsets of $G$, and there are exponentially many of them.
Therefore the following result of Thomason, providing a sufficient
condition for pseudo-randomness based on degrees and co-degrees only,
carries a certain element of surprise in it.
\begin{theo}\label{jumlocal}\cite{Tho87a}
Let $G$ be a graph on $n$ vertices with minimum degree $np$. If no pair
of vertices of $G$ has more than $np^2+l$ common neighbors, then $G$ is
$(p,\sqrt{(p+l)n})$-jumbled. 
\end{theo}
The above theorem shows how the pseudo-randomness condition of
(\ref{jumbled}) can be ensured/checked by testing only a polynomial
number of easily accessible conditions. It is very useful for showing
that specific constructions are jumbled. Also, it can find algorithmic
applications, for example, a very similar approach has been used by 
Alon, Duke, Lefmann, R\"odl and Yuster in their Algorithmic Regularity
Lemma \cite{AloDukLefRodYus94}. 

As observed by Thomason, the minimum degree condition of Theorem
\ref{jumlocal} can be dropped if we require that every pair of vertices
has $(1+o(1))np^2$ common neighbors. One cannot however weaken the
conditions of the
theorem so as to only require that every {\em edge} is in at most
$np^2+l$ triangles.

Another sufficient condition for pseudo-randomness, this time of global
nature, has also been provided in \cite{Tho87a}, \cite{Tho87b}:
\begin{theo}\label{jumglobal}\cite{Tho87a}
Let $G$ be a graph of order $n$, let $\eta n$ be an integer between 2
and $n-2$, and let $\omega>1$ be a real number. Suppose that each
induced subgraph $H$ of order $\eta n$ satisfies 
$|e(H)-p{{\eta n}\choose 2}|\le \eta n\alpha$. Then $G$ is
$(p,7\sqrt{n\alpha/\eta}/(1-\eta))$-jumbled. Moreover $G$ contains a
subset $U\subseteq V(G)$ of size 
$|U|\ge \left(1-\frac{380}{n(1-\eta)^2w}\right)n$ such that the induced
subgraph $G[U]$ is $(p,\omega\alpha)$-jumbled.
\end{theo}

Thomason also describes in \cite{Tho87a}, \cite{Tho87b} several
properties of jumbled graphs. We will not discuss these results in
details here as we will mostly adopt a different approach to
pseudo-randomness. Occasionally however we will compare some of later
results to those obtained by Thomason.  

\subsection{Equivalent definitions of weak pseudo-randomness}
Let us go back to the jumbledness condition (\ref{jumbled}) of Thomason.
As we have already noted it becomes non-trivial only when the error
term in (\ref{jumbled}) is $o(n^2p)$. Thus the latter condition can be
considered as the weakest possible condition for pseudo-randomness.

Guided by the above observation we now define the notion of weak
pseudo-randomness as follows. Let $(G_n)$ be a sequence of graphs,
where $G_n$ has $n$ vertices. Let also $p=p(n)$ is a parameter ($p(n)$
is a typical density of graphs in the sequence). We say that the
sequence $(G_n)$ is {\em weakly pseudo-random} if the following
condition holds:
\begin{eqnarray}\label{weakpr}
\mbox{For all subsets $U\subseteq V(G_n)$,}\quad\quad 
     \left|e(U)-p{{|U|}\choose 2}\right|=o(n^2p)\ .
\end{eqnarray}
For notational convenience we will frequently write $G=G_n$, tacitly
assuming that $(G)$ is in fact a sequence of graphs.

Notice that the error term in the above condition of weak
pseudo-randomness does not depend on the size of the subset $U$.
Therefore it applies essentially only to subsets $U$ of linear size,
ignoring subsets $U$ of size $o(n)$. Hence (\ref{weakpr}) is potentially
much weaker than Thomason's jumbledness condition (\ref{jumbled}).  

Corollary \ref{dist3} supplies us with the first example of weakly
pseudo-random graphs -- a random graph $G(n,p)$ is weakly pseudo-random
as long as $p(n)$ satisfies $np\rightarrow\infty$. We can thus say that
if a graph $G$ on $n$ vertices is weakly pseudo-random for a parameter 
$p$, then the edge distribution of $G$ is close to that of $G(n,p)$.

In the previous subsection we have already seen examples of conditions
implying pseudo-randomness. In general one can expect that conditions of
various kinds that hold almost surely in $G(n,p)$ may imply or be
equivalent to weak pseudo-randomness of graphs with edge density $p$.

Let us first consider the case of the constant edge density $p$. This 
case has been treated extensively in the celebrated paper of Chung,
Graham and Wilson from 1989 \cite{ChuGraWil89}, where they formulated
several equivalent conditions for weak pseudo-randomness. In order to
state their important result we need to introduce some notation.

Let $G=(V,E)$ be a graph on  $n$ vertices.
For a graph $L$ we denote by $N^*_G(L)$ the
number of labeled induced copies of $L$ in $G$, and by $N_G(L)$ the
number of labeled not necessarily induced copies of $L$ in $G$. For a
pair of vertices $x,y\in V(G)$, we set $s(x,y)$ to be the number of
vertices of $G$ joined to $x$ and $y$ the same way: either to both or to
none. Also, $codeg(x,y)$ is the number of common neighbors of $x$ and
$y$ in $G$. Finally, we order the eigenvalues $\lambda_i$ of the
adjacency matrix $A(G)$ so that $|\lambda_1|\ge |\lambda_2|\ge\ldots\ge
|\lambda_n|$.

\begin{theo}\label{CGW}\cite{ChuGraWil89}
Let $p\in (0,1)$ be fixed. For any graph sequence $(G_n)$ the following
properties are equivalent:
\begin{description}
\item[$P_1(l)$:\quad] 
For a fixed $l\ge 4$ for all graphs $L$ on $l$ vertices,
$$
N_G^*(L)=(1+o(1))n^lp^{|E(L)|}(1-p)^{{l\choose 2}-|E(L)|}\ .
$$
\item[$P_2(t)$:\quad] 
Let $C_t$ denote the cycle of length $t$. Let $t\ge 4$ be even,
$$
e(G_n)=\frac{n^2p}{2}+o(n^2)\quad\mbox{and}\quad 
N_G(C_t)\le (np)^t+o(n^t)\ .
$$
\item[$P_3$:\quad] 
$
e(G_n)\ge \frac{n^2p}{2}+o(n^2)\quad\mbox{and}\quad
\lambda_1=(1+o(1))np,~~ \lambda_2=o(n)\ .
$
\item[$P_4$:\quad] For each subset $U\subset V(G)$,\quad
$e(U)=\frac{p}{2}|U|^2+o(n^2)$\ .    
\item[$P_5$:\quad] For each subset $U\subset V(G)$ with 
$|U|=\lfloor\frac{n}{2}\rfloor$,\quad we have\quad 
$e(U)=\left(\frac{p}{8}+o(1)\right)n^2\ .$
\item[$P_6$:\quad]
$\sum_{x,y\in V} |s(x,y)-(p^2+(1-p)^2)n|=o(n^3)$\ .
\item[$P_7$:\quad]
$\sum_{x,y\in V} |codeg(x,y)-p^2n|=o(n^3)$\ .
\end{description}
\end{theo}

Note that condition $P_4$ of this remarkable theorem is in fact
identical to our condition (\ref{weakpr}) of weak pseudo-randomness.
Thus according to the theorem all conditions $P_1$--$P_3$, $P_5-P_7$ are
in fact equivalent to weak pseudo-randomness!

As noted by Chung et al. probably the most surprising fact
(although possibly less surprising for the reader in view of Theorem  
\ref{jumlocal}) is that apparently the weak condition $P_2(4)$ is strong
enough to imply weak pseudo-randomness.

It is quite easy to add another condition to the equivalence list of the
above theorem: for all $U,W\subset V$, $e(U,W)=p|U||W|+o(n^2)$. 

A condition of a very different type, related to the celebrated
Szemer\'edi Regularity Lemma  has been added to the above list by
Simonovits and S\'os in \cite{SimSos91}. They showed that if a graph $G$
possesses a Szemer\'edi partition in which almost all pairs have density
$p$, then $G$ is weakly pseudo-random, and conversely if $G$ is weakly
pseudo-random then in every Szemer\'edi partition all pairs are regular
with density $p$. An extensive background on the Szemer\'edi Regularity
Lemma, containing in particular the definitions of the above used
notions, can be found in a survey paper of Koml\'os and Simonovits
\cite{KomSim96}. 

The reader may have gotten the feeling that basically every property of
random graphs $G(n,p)$ ensures weak pseudo-randomness. This feeling is
quite misleading, and one should be careful while formulating properties
equivalent to pseudo-randomness. Here is an example provided by Chung et
al. Let $G$ be a graph with vertex set $\{1,\ldots,4n\}$ defined as
follows: the subgraph of $G$ spanned by the first $2n$ vertices is a
complete bipartite graph $K_{n,n}$, the subgraph spanned by the last
$2n$ vertices is the complement of $K_{n,n}$, and for every pair
$(i,j),1\le i\le 2n, 2n+1\le j\le 4n$, the edge $(i,j)$ is present in
$G$ independently with probability $0.5$. Then $G$ is almost surely
a graph on $4n$
vertices with edge density $0.5$. One can verify that $G$ has properties
$P_1(3)$ and $P_2(2t+1)$ for every $t\ge 1$, but is obviously very far
from being pseudo-random (contains a clique and an independent set of
one quarter of its size). Hence $P_1(3)$ and $P_2(2t+1)$ are not
pseudo-random properties. This example shows also the real difference
between even and odd cycles in this context -- recall that Property
$P_2(2t)$ does imply pseudo-randomness.  

A possible explanation to the above described somewhat disturbing
phenomenon has been suggested by Simonovits and S\'os in
\cite{SimSos97}. They noticed that the above discussed properties are
not hereditary in the sense that the fact that the whole graph $G$
possesses one of these properties does not imply that large induced
subgraphs of $G$ also have it. A property is called {\em hereditary} in
this context if it is assumed to hold for all sufficiently large
subgraphs $F$ of our graph $G$ with the same error term as for $G$.
Simonovits and S\'os proved that adding this hereditary condition gives
significant extra strength to many properties making them
pseudo-random.
\begin{theo}\cite{SimSos97}\label{SSher}
Let $L$ be a fixed graph on $l$ vertices, and let $p\in (0,1)$ be fixed.
Let $(G_n)$ be a sequence of graphs. If for every induced subgraph
$H\subseteq G$ on $h$ vertices,
$$
N_H(L)=p^{|E(L)|}h^l+o(n^l)\,,
$$ 
then $(G_n)$ is weakly pseudo-random, i.e. property $P_4$ holds.
\end{theo}
Two main distinctive features of the last result compared to Theorem
\ref{CGW} are: (a) $P_1(3)$ assumed hereditarily implies
pseudo-randomness; and (b) requiring the right number of copies of a
{\em single} graph $L$ on $l$ vertices is enough, compared to Condition
$P_1(l)$ required to hold for {\em all} graphs on $l$ vertices
simultaneously.

Let us switch now to the case of vanishing edge density $p(n)=o(1)$.
This case has been treated in two very recent papers of Chung and Graham
\cite{ChuGra02} and of Kohayakawa, R\"odl and Sissokho
\cite{KohRodSis02}.
Here the picture becomes significantly more complicated compared to the
dense case. In particular, there exist graphs with very balanced edge
distribution not containing a single copy of some fixed subgraphs (see
the Erd\H os-R\'enyi graph and the Alon graph in the next section
(Examples 6, 9, resp.)). 

In an attempt to find properties equivalent to weak pseudo-randomness in
the sparse case, Chung and Graham define the following properties in
\cite{ChuGra02} :

\noindent{\bf CIRCUIT($t$):} The number of closed walks
$w_0,w_1,\ldots,w_t=w_0$ of length $t$ in $G$ is $(1+o(1))(np)^t$;

\noindent{\bf CYCLE($t$):} The number of labeled $t$-cycles in $G$ is
$(1+o(1))(np)^t$;

\noindent{\bf EIG:} The eigenvalues $\lambda_i$,
$|\lambda_1|\ge|\lambda_2|\ge\ldots |\lambda_n|$, of the adjacency
matrix of $G$ satisfy:
\begin{eqnarray*}
\lambda_1&=&(1+o(1))np\,,\\
|\lambda_i|&=&o(np), i>1\,.
\end{eqnarray*}

\noindent{\bf DISC:} For all $X,Y\subset V(G)$,
$$
|e(X,Y)-p|X||Y||=o(pn^2)\ .
$$

(DISC here is in fact DICS(1) in \cite{ChuGra02}). 

\begin{theo}\label{CG1}\cite{ChuGra02}
Let $(G=G_n: n\rightarrow \infty)$ be a sequence of graphs with
$e(G_n)=(1+o(1))p{n\choose 2}$. Then the following implications hold for
all $t\ge 1$:
$$
CIRCUIT(2t)\Rightarrow EIG\Rightarrow DISC\ .
$$
\end{theo}
 
\noindent{\bf Proof.\quad} 
To prove the first implication, let $A$ be the adjacency matrix of $G$,
and consider the trace $Tr(A^{2t})$. The $(i,i)$-entry of $A^{2t}$ is
equal to the number of closed walks of length $2t$
starting and ending at $i$, and hence $Tr(A^{2t})=(1+o(1))(np)^{2t}$.
On the other hand, since $A$ is symmetric it is similar to the diagonal
matrix $D=diag(\lambda_1,\lambda_2,\ldots,\lambda_n)$, and therefore 
$Tr(A^{2t})=\sum_{i=1}^{2t}\lambda_i^{2t}$. We obtain:
$$
\sum_{i=1}^n\lambda_i^{2t}=(1+o(1))(np)^{2t}\ .
$$
Since the first eigenvalue of $G$ is easily shown to be as large as its
average degree, it follows that $\lambda_1\ge
2|E(G)|/|V(G)|=(1+o(1))np$. Combining these two facts we derive that
$\lambda_1=(1+o(1))np$ and $|\lambda_i|=o(np)$ as required. 

The second implication will be proven in the next subsection. \hfill
$\Box$

\medskip

Both reverse implications are false in general. To see why
$DISC\not\Rightarrow EIG$ take a graph $G_0$ on $n-1$ vertices with all
degrees equal to $(1+o(1))n^{0.1}$ and having property $DISC$ (see next
section for examples of such graphs). Now add to $G_0$ a vertex $v^*$
and connect it to any set of size $n^{0.8}$ in $G_0$, let $G$
be the obtained graph. Since $G$ is obtained from $G_0$ by adding
$o(|E(G_0|)$ edges, $G$ still satisfies $DISC$. On the other hand, $G$
contains a star $S$ of size $n^{0.8}$ with a center at $v^*$,
and hence $\lambda_1(G)\ge\lambda_1(S)=\sqrt{n^{0.8}-1}\gg |E(G)|/n$
(see, e.g. Chapter 11 of \cite{Lov93} for the relevant proofs). This
solves an open question from \cite{ChuGra02}.  

The Erd\H os-R\'enyi graph from the next section is easily seen to
satisfy $EIG$, but fails to satisfy $CIRCUIT(4)$. Chung and Graham
provide an alternative example in \cite{ChuGra02} (Example 1).

The above discussion indicates that one probably needs to impose some
additional condition on the graph $G$ to glue all these pieces together
and to make the above stated properties equivalent. One such condition
has been suggested by Chung and Graham who defined:

\medskip
\noindent{\bf U($t$):} For some absolute constant $c$, all degrees in
$G$ satisfy: $d(v)< cnp$, and for every pair of vertices $x,y\in G$
the number
$e_{t-1}(x,y)$ of walks of length $t-1$ from $x$ to $y$ satisfies:
$e_{t-1}(x,y)\le cn^{t-2}p^{t-1}$.
\medskip

Notice that $U(t)$ can only hold for $p>c'n^{-1+1/(t-1)}$, where $c'$
depends on $c$. Also, every dense graph ($p=\Theta(1)$) satisfies
$U(t)$. 

As it turns out adding property $U(t)$ makes all the above defined
properties equivalent and thus equivalent to the notion of weak
pseudo-randomness (that can be identified with property $DISC$): 

\begin{theo}\label{CG2}\cite{ChuGra02}
Suppose for some constant $c>0$, $p(n)>cn^{-1+1/(t-1)}$, where $t\ge 2$.
For any family of graphs $G_n$, $|E(G_n)|=(1+o(1))p{n\choose 2}$,
satisfying $U(t)$, the following properties are all equivalent:
$CIRCUIT(2t), CYCLE(2t), EIG$ and $DISC$.
\end{theo}

Theorem \ref{CG2} can be viewed as a sparse analog of Theorem \ref{CGW}
as it also provides a list of conditions equivalent to weak
pseudo-randomness.

Further properties implying or equivalent to pseudo-randomness,
including local statistics conditions, are given in \cite{KohRodSis02}.

\subsection{Eigenvalues and pseudo-random graphs}
In this subsection we describe an approach to pseudo-randomness based on
graph eigenvalues -- the approach most frequently used in this survey.
Although the eigenvalue-based condition
is not as general as the jumbledness
condition of Thomason or some other properties described in the previous
subsection, its power and convenience are so appealing that they
certainly constitute a good enough reason to prefer this approach. Below
we first provide a necessary background on graph spectra
and then derive quantitative estimates connecting the eigenvalue gap
and edge distribution. 

Recall that the {\em adjacency matrix} of a graph $G=(V,E)$ with vertex
set $V=\{1,\ldots,n\}$ is an $n$-by-$n$ matrix whose entry $a_{ij}$ is
1 if $(i,j)\in E(G)$, and is 0 otherwise. Thus $A$ is a $0,1$ symmetric
matrix with zeroes along the main diagonal, and we can apply the
standard machinery of eigenvalues and eigenvectors of real symmetric
matrices. It follows that all eigenvalues of $A$ (usually also
called the eigenvalues of the graph $G$ itself) are real, and we denote
them by $\lambda_1\ge\lambda_2\ge\ldots\ge \lambda_n$. Also, there is
an orthonormal basis $B=\{x_1,\ldots,x_n\}$ of the euclidean space $R^n$
composed of eigenvectors of $A$: $Ax_i=\lambda_i x_i$, $x_i^tx_i=1$, 
$i=1,\ldots,n$. The matrix $A$ can be decomposed then as:
$A=\sum_{i=1}^n\lambda_ix_ix_i^t$ -- the so called spectral
decomposition of $A$. (Notice that the product $xx^t$, $x\in R^n$, is an
$n$-by-$n$ matrix of rank 1; if $x,y,z\in R^n$ then
$y^t(xx^t)z=(y^tx)(x^tz)$). Every vector $y\in R^n$ can be easily
represented in basis $B$: $y=\sum_{i=1}^n(y^tx_i)x_i$. Therefore, for
$y,z\in R^n$, $y^tz=\sum_{i=1}^n (y^tx_i)(z^tx_i)$ and
$\|y\|^2=y^ty=\sum_{i=1}^n (y^tx_i)^2$.  

All the above applies in fact to all real symmetric matrices.
Since the adjacency matrix $A$ of a graph $G$ is a matrix with
non-negative entries, one can derive some important extra features
of $A$, most notably the Perron-Frobenius Theorem, that reads in
the graph context as follows: if $G$ is connected then the
multiplicity of $\lambda_1$ is one, all coordinates of the first
eigenvector $x_1$ can be assumed to be strictly positive, and
$|\lambda_i|\le \lambda_1$ for all $i\ge 2$. Thus, graph spectrum
lies entirely in the interval $[-\lambda_1,\lambda_1]$. 

For the most important special case of regular graphs
Perron-Frobenius implies the following corollary:

\begin{prop}\label{PF}
Let $G$ be a $d$-regular graph on $n$ vertices. Let
$\lambda_1\ge\lambda_2\ge\ldots\ge \lambda_n$ be the eigenvalues
of $G$. Then $\lambda_1=d$ and $-d\le \lambda_i \le d$ for all
$1\le i\le n$. Moreover, if $G$ is connected then the first
eigenvector $x_1$ is proportional to the all one vector 
$(1,\ldots,1)^t\in R^n$, and $\lambda_i<d$ for all $i\ge 2$.
\end{prop}

To derive the above claim from the Perron-Frobenius Theorem
observe that $e=(1,\ldots,1)$ is immediately seen to be an
eigenvector of $A(G)$ corresponding to the eigenvalue $d$: $Ae=de$.
The positivity of the coordinates of $e$ implies then that $e$ is not 
orthogonal to the first eigenvector, and hence is
in fact proportional to $x_1$ of $A(G)$.
Proposition \ref{PF} can be also proved directly without relying
on the Perron-Frobenius Theorem.

We remark that $\lambda_n=-d$ is possible, in fact it holds if and
only if the graph $G$ is bipartite.

All this background information, presented above in a somewhat
condensed form, can be found in many textbooks in Linear  Algebra.
Readers more inclined to consult combinatorial books can find it
for example in a recent monograph of Godsil and Royle on Algebraic
Graph Theory \cite{GodRoy01}.

We now prove a well known theorem (see its variant, e.g., in Chapter 9,
\cite{AloSpe00}) bridging between graph spectra and edge distribution.

\begin{theo}\label{eigen}
Let $G$ be a $d$-regular graph on $n$ vertices. Let
$d=\lambda_1\ge\lambda_2\ge\ldots\lambda_n$ be the eigenvalues of $G$.
Denote 
$$
\lambda=max_{2\le i\le n}|\lambda_i|\,.
$$ 
Then for every two subsets $U,W\subset V$,
\begin{equation}\label{eig}
\left|e(U,W)-\frac{d|U||W|}{n}\right|\le 
\lambda\sqrt{|U||W|\left(1-\frac{|U|}{n}\right)
\left(1-\frac{|W|}{n}\right)}\ .
\end{equation}
\end{theo}

\noindent{\bf Proof.\ } Let $B=\{x_1,\ldots,x_n\}$ be an orthonormal
basis of $R^n$ composed from eigenvectors of $A$: $Ax_i=\lambda_ix_i$,
$1\le i\le n$. We represent $A=\sum_{i=1}^n\lambda_ix_ix_i^t$. Denote 
\begin{eqnarray*}
A_1&=&\lambda_1x_1x_1^t\,,\\
\Epsilon&=&\sum_{i=2}^n\lambda_ix_ix_i^t\,, 
\end{eqnarray*}
then $A=A_1+\Epsilon$.

Let $u=|U|$, $w=|W|$ be the cardinalities of $U,W$, respectively. We
denote the characteristic vector of $U$ by $\chi_U\in R^n$, i.e.
$\chi_U(i)=1$ if $i\in U$, and $\chi_U(i)=0$ otherwise. Similarly, let
$\chi_W\in R^n$ be the characteristic vector of $W$. We represent
$\chi_U$, $\chi_W$ according to $B$:
\begin{eqnarray*}
\chi_U &=& \sum_{i=1}^n \alpha_ix_i,\quad  \alpha_i=\chi_U^tx_i,
          \quad \sum_{i=1}^n\alpha_i^2=\|\chi_U\|^2=u\,,\\
\chi_W &=& \sum_{i=1}^n \beta_ix_i,  \quad \beta_i=\chi_W^tx_i,
          \quad \sum_{i=1}^n\beta_i^2=\|\chi_W\|^2=w\ .
\end{eqnarray*}
It follows easily from the definitions of $A$, $\chi_U$ and $\chi_W$
that the product $\chi_U^tA\chi_W$ counts exactly the number of edges of
$G$ with one endpoint in $U$ and the other one in $W$, i.e. 
$$
e(U,W)=\chi_U^tA\chi_W\ =\chi_U^tA_1\chi_W+\chi_U^t\Epsilon\chi_W\ .
$$
Now we estimate the last two summands separately, the first of them
will be the main term for $e(U,W)$, the second one will be the error
term. Substituting the expressions for $\chi_U$, $\chi_W$ and recalling
the orthonormality of $B$, we get:
\begin{equation}\label{fir}
\chi_U^tA_1\chi_W=\left(\sum_{i=1}^n\alpha_ix_i\right)^t
                  (\lambda_1x_1x_1^t)
                  \left(\sum_{j=1}^n\beta_jx_j\right)=
                  \sum_{i=1}^n\sum_{j=1}^n \alpha_i\lambda_1\beta_j
                  (x_i^tx_1)(x_1^tx_j)=\alpha_1\beta_1\lambda_1\ .
\end{equation}
Similarly,
\begin{equation}\label{rest}
\chi_U^t\Epsilon\chi_W=\left(\sum_{i=1}^n\alpha_ix_i\right)^t
\left(\sum_{j=2}^n\lambda_jx_jx_j^t\right)
\left(\sum_{k=1}^n\beta_kx_k\right)=
\sum_{i=2}^n\alpha_i\beta_i\lambda_i\ .
\end{equation}

Recall now that $G$ is $d$-regular. Then according to Proposition
\ref{PF}, $\lambda_1=d$ and $x_1=\frac{1}{\sqrt{n}}(1,\ldots,1)^t$. We
thus get: $\alpha_1=\chi_U^tx_1=u/\sqrt{n}$ and
$\beta_1=\chi_W^tx_1=w/\sqrt{n}$. Hence it follows from (\ref{fir}) that
$\chi_U^tA_1\chi_W=duw/n$.

Now we estimate the absolute value of the error term 
$\chi_U^t\Epsilon\chi_W$. Recalling (\ref{rest}), the definition of
$\lambda$ and the obtained values of $\alpha_1$, $\beta_1$, we
derive, applying Cauchy-Schwartz:
\begin{eqnarray*}
|\chi_U^t\Epsilon\chi_W|&=&|\sum_{i=2}^n\alpha_i\beta_i\lambda_i|
\le\lambda|\sum_{i=2}^n\alpha_i\beta_i|\le 
\lambda\sqrt{\sum_{i=2}^n\alpha_i^2\sum_{i=2}^n\beta_i^2}\\
&=&\lambda\sqrt{(\|\chi_U\|^2-\alpha_1^2)(\|\chi_W\|^2-\beta_1^2)}=
\lambda\sqrt{\left(u-\frac{u^2}{n}\right)\left(w-\frac{w^2}{n}\right)}
\ .
\end{eqnarray*}
The theorem follows.\hfill $\Box$

\bigskip

The above proof can be extended to the irregular (general) case.
Since the obtained quantitative bounds on edge distribution turn out to
be somewhat cumbersome, we will just indicate how they can be obtained.
Let $G=(V,E)$ be a graph on $n$ vertices with {\em average} degree $d$.
Assume that the eigenvalues of $G$ satisfy $\lambda<d$, with $\lambda$
as defined in the theorem. Denote 
$$
K=\sum_{v\in V}(d(v)-d)^2\ .
$$
The parameter $K$ is a measure of irregularity of $G$. Clearly $K=0$
if and only if $G$ is $d$-regular. Let
$e=\frac{1}{\sqrt{n}}(1,\ldots,1)^t$. We represent $e$ in the basis
$B=\{x_1,\ldots,x_n\}$ of the eigenvectors of $A(G)$:
$$
e=\sum_{i=1}^n\gamma_ix_i,\quad \gamma_i=e^tx_i,\quad 
\sum_{i=1}^n\gamma_i^2=\|e\|^2=1\ .
$$
Denote $z=\frac{1}{\sqrt{n}}(d(v_1)-d,\ldots,d(v_n)-d)^t$, then
$\|z\|^2=K/n$. Notice that
$Ae=\frac{1}{\sqrt{n}}(d(v_1),\ldots,d(v_n))^t$ $=de+z$, and therefore 
$z=Ae-de=\sum_{i=1}^n\gamma_i(\lambda_i-d)x_i$. This implies:
\begin{eqnarray*}
\frac{K}{n}&=&\|z\|^2= \sum_{i=1}^n\gamma_i^2(\lambda_i-d)^2\ge 
\sum_{i=2}^n\gamma_i^2(\lambda_i-d)^2\\
&\ge& (d-\lambda)^2\sum_{i=2}^n\gamma_i^2\ .
\end{eqnarray*}
Hence $\sum_{i=2}^n\gamma_i^2\le \frac{K}{n(d-\lambda)^2}$. It follows
that $\gamma_1^2=1-\sum_{i=2}^n\gamma_i^2\ge 1-\frac{K}{n(d-\lambda)^2}$
and 
$$
\gamma_1\ge \gamma_1^2 \ge 1-\frac{K}{n(d-\lambda)^2}\ .
$$
Now we estimate the distance between the vectors $e$ and $x_1$ and show
that they are close given that the parameter $K$ is small. 
\begin{eqnarray*}
\|e-x_1\|^2&=& (e-x_1)^t(e-x_1)=e^te+x_1^tx_1-2e^tx_1=1+1-2\gamma_1
            =2-2\gamma_1\\
           &\le& \frac{2K}{n(d-\lambda)^2}\ .
\end{eqnarray*}

We now return to expressions (\ref{fir}) and (\ref{rest}) from the proof
of Theorem \ref{eigen}. In order to estimate the main term
$\chi_U^tA_1\chi_W$, we bound the coefficients $\alpha_1$, $\beta_1$
and $\lambda_1$ as follows:
$$
\alpha_1=\chi_U^tx_1=\chi_U^te+\chi_U^t(x_1-e)= 
         \frac{u}{\sqrt{n}}+\chi_U^t(x_1-e)\ ,
$$
and therefore 
\begin{equation}\label{alpha1}
\left|\alpha_1-\frac{u}{\sqrt{n}}\right|=|\chi_U^t(x_1-e)|
\le \|\chi_U||\cdot\|x_1-e\|\le \frac{\sqrt{\frac{2Ku}{n}}}{d-\lambda} 
\ .
\end{equation}
In a similar way one gets:
\begin{equation}\label{beta1}
\left|\beta_1-\frac{w}{\sqrt{n}}\right|\le 
                                \frac{\sqrt{\frac{2Kw}{n}}}{d-\lambda}
\ .
\end{equation}
Finally, to estimate from above the absolute value of the difference
between $\lambda_1$ and $d$ we argue as follows:
$$
\frac{K}{n}=\|z\|^2=\sum_{i=1}^n\gamma_i^2(\lambda_i-d)^2\ge
\gamma_1^2(\lambda_1-d)^2\,,
$$
and therefore 
\begin{equation}\label{lambda1}
|\lambda_1-d|\le \frac{1}{\gamma_1}\sqrt{\frac{K}{n}}
             \le \frac{n(d-\lambda)^2}{n(d-\lambda)^2-K}
                  \sqrt{\frac{K}{n}}\ .
\end{equation}
Summarizing, we see from (\ref{alpha1}), (\ref{beta1}) and
(\ref{lambda1}) that the main term in the product $\chi_U^tA_1\chi_W$ is
equal to $\frac{duw}{n}$, just as in the regular case, and the error
term is governed by the parameter $K$.

In order to estimate the error term $\chi_U^t\Epsilon \chi_W$ we use
(\ref{rest}) to get:
\begin{eqnarray*}
\hspace{3.1cm}
|\chi_U^t\Epsilon\chi_W|&=&
\left|\sum_{i=2}^n\alpha_i\beta_i\lambda_i\right|
\le \lambda\left|\sum_{i=2}^n\alpha_i\beta_i\right|
\le \lambda \sqrt{\sum_{i=2}^n\alpha_i^2\sum_{i=2}^n\beta_i^2}\\
&\le&   \lambda \sqrt{\sum_{i=1}^n\alpha_i^2\sum_{i=1}^n\beta_i^2}
= \lambda \|\chi_U\|\,\|\chi_W\|=\lambda\sqrt{uw}. \hspace{3.1cm} \Box
\end{eqnarray*}

Applying the above developed techniques we can prove now the second
implication of Theorem \ref{CG1}. Let us prove first that $EIG$ implies
$K=o(nd^2)$, where $d=(1+o(1))np$ is as before the average degree of
$G$. Indeed, for every vector $v\in R^n$ we have $\|Av\|\le
\lambda_1\|v\|$, and therefore 
$$
\lambda_1^2n=\lambda_1^2e^te\ge (Ae)^t(Ae)=\sum_{v\in V}d^2(v)\ .
$$
Hence from $EIG$ we get: $\sum_{v\in V}d^2(v)\le (1+o(1))nd^2$. 
As $\sum_{v}d(v)=nd$, it follows that:
$$
K=\sum_{v\in V}(d(v)-d)^2=\sum_{v\in V}d^2(v)-2d\sum_{v\in V}d(v)+nd^2
 = (1+o(1))nd^2-2nd^2+nd^2=o(nd^2)\,,
$$
 as promised. Substituting this into
estimates (\ref{alpha1}), (\ref{beta1}), (\ref{lambda1}) and using
$\lambda=o(d)$ of $EIG$ we get:
\begin{eqnarray*}
\alpha_1  &=& \frac{u}{\sqrt{n}}+o(\sqrt{u})\,,\\
\beta_1  &=&  \frac{w}{\sqrt{n}}+o(\sqrt{w})\,,\\
\lambda_1 &=& (1+o(1))d\,,
\end{eqnarray*}
and therefore 
$$
\chi_U^tA_1\chi_W=\frac{duw}{n}+o(dn)\ .
$$
Also, according to $EIG$, $\lambda=o(d)$, which implies:
$$
\chi_U^t\Epsilon\chi_w=o(d\sqrt{uw})=o(dn)\,,
$$
and the claim follows. \hfill $\Box$
 
Theorem \ref{eigen} is a truly remarkable result. Not only it connects
between two seemingly unrelated graph characteristics -- edge
distribution and spectrum, it also provides a very good quantitative
handle for the uniformity of edge distribution, based on easily
computable, both theoretically and practically, graph parameters --
graph eigenvalues. According to the bound (\ref{eig}), a polynomial number
of parameters can control quite well the number of edges in 
exponentially many subsets of vertices.

The parameter $\lambda$ in the formulation of Theorem \ref{eigen} is
usually called the {\em second eigenvalue} of the $d$-regular graph $G$
(the first and the trivial one being $\lambda_1=d$). There is
certain inaccuracy though in this term, as in fact 
$\lambda=\max\{\lambda_2,-\lambda_n\}$. Later we will call, following 
Alon, a $d$-regular
graph $G$ on $n$ vertices in which all eigenvalues, but the first one,
are at most $\lambda$ in their absolute values, an {\em 
$(n,d,\lambda)$-graph}. 

Comparing (\ref{eig}) with the definition of jumbled graphs by Thomason
we see that an $(n,d,\lambda)$-graph $G$ is $(d/n,\lambda)$-jumbled. Hence
the parameter $\lambda$ (or in other words, the so called {\em spectral
gap} -- the difference between $d$ and $\lambda$) is responsible for
pseudo-random properties of such a graph. The smaller the value of
$\lambda$ compared to $d$, the more close is the edge distribution of
$G$ to the ideal uniform distribution. A natural question is then: how
small can be $\lambda$? It is easy to see that as long as $d\le
(1-\epsilon)n$, $\lambda=\Omega(\sqrt{d})$. Indeed, the trace of $A^2$
satisfies:
$$
nd=2|E(G)|=Tr(A^2)=\sum_{i=1}^n\lambda_i^2\le d^2+(n-1)\lambda_2
\le (1-\epsilon)nd+(n-1)\lambda^2\,,
$$
and $\lambda=\Omega(\sqrt{d})$ as claimed. More accurate bounds are
known for smaller values of $d$ (see, e.g. \cite{Nil91}). Based on these
estimates we can say that an $(n,d,\lambda)$-graph $G$, for which
$\lambda=\Theta(\sqrt{d})$, is a very good pseudo-random graph. We will
see several examples of such graphs in the next section.

\subsection{Strongly regular graphs}
A {\em strongly regular graph} $srg(n,d,\eta,\mu)$ is a $d$-regular
graph on $n$ vertices in which every pair of adjacent vertices has
exactly $\eta$ common neighbors and every pair of non-adjacent vertices
has exactly $\mu$ common neighbors. (We changed the very standard notation
in the above definition so as to avoid interference with other
notational conventions throughout this paper and to make it more
coherent,
usually the parameters are denoted $(v,k,\lambda,\mu)$). Two simple
examples of strongly regular graph are the pentagon $C_5$ that has
parameters $(5,2,0,1)$, and the Petersen graph whose parameters are
$(10,3,0,1)$. Strongly regular graphs were introduced by Bose in 1963
\cite{Bos63} who also pointed out their tight connections with finite
geometries. As follows from the definition, strongly regular graphs
are highly
regular structures, and one can safely predict that algebraic methods
are extremely useful in their study. We do not intend to provide any
systematic coverage of this fascinating concept here, addressing the
reader to the vast literature on the subject instead (see, e.g.,
\cite{BrovLi84}). Our aim here is to calculate the eigenvalues of
strongly regular graphs and then to connect them with pseudo-randomness,
relying on results from the previous subsection.

\begin{prop}\label{srg}
Let $G$ be a connected strongly regular graph with parameters 
$(n,d,\eta,\mu)$. Then the eigenvalues of $G$ are: $\lambda_1=d$ with
multiplicity $s_1=1$,
$$
\lambda_2=\frac{1}{2}\left(\eta-\mu+\sqrt{(\eta-\mu)^2+4(d-\mu)}\right) 
$$ 
and 
$$
\lambda_3=\frac{1}{2}\left(\eta-\mu-\sqrt{(\eta-\mu)^2+4(d-\mu)}\right)
\,,
$$
with multiplicities 
$$
s_2=\frac{1}{2}\left(n-1+\frac{(n-1)(\mu-\eta)-2d}
                              {\sqrt{(\mu-\eta)^2+4(d-\mu)}}\right)
$$
and 
$$
s_3=\frac{1}{2}\left(n-1-\frac{(n-1)(\mu-\eta)-2d}
                              {\sqrt{(\mu-\eta)^2+4(d-\mu)}}\right)\,,
$$
respectively.
\end{prop}

\noindent{\bf Proof.\ } Let $A$ be the adjacency matrix of $A$. By
the definition of $A$ and the fact that $A$ is symmetric with zeroes on
the main diagonal, the $(i,j)$-entry of the square $A^2$ counts the
number of common neighbors of $v_i$ and $v_j$ in $G$ if $i\ne j$, and is
equal to the degree $d(v_i)$ in case $i=j$. 
The statement that $G$ is $srg(n,d,\eta,\mu)$ is equivalent then to:
\begin{equation}\label{srg1}
AJ=dJ,\quad\quad A^2=(d-\mu)I+\mu J+(\eta-\mu)A\ ,
\end{equation}
where $J$ is the $n$-by-$n$ all-one matrix and $I$ is the $n$-by-$n$
identity matrix.

Since $G$ is $d$-regular and connected, we obtain from the
Perron-Frobenius Theorem that $\lambda_1=d$ is an eigenvalue of $G$ with
multiplicity 1 and with $e=(1,\ldots,1)^t$ as the corresponding
eigenvector. Let $\lambda\ne d$ be another eigenvalue of $G$, and let
$x\in R^n$ be a corresponding eigenvector. Then $x$ is orthogonal to
$e$, and therefore $Jx=0$. Applying both sides of the second identity in
(\ref{srg1}) to $x$ we get the equation:
$\lambda^2x=(d-\mu)x+(\eta-\mu)\lambda x$, which results in the
following quadratic equation for $\lambda$:
$$
\lambda^2+(\mu-\eta)\lambda+(\mu-d)=0\ .
$$
This equation has two solutions $\lambda_2$ and $\lambda_3$ as defined
in the proposition formulation. If we denote by $s_2$ and $s_3$ the
respective multiplicities of $\lambda_2$ and $\lambda_3$ as eigenvalues
of $A$, we get:
$$
1+s_2+s_3=n,\quad\quad Tr(A)=d+s_2\lambda_2+s_3\lambda_3=0\ .
$$
Solving the above system of linear equations for $s_2$ and $s_3$ we obtain
the assertion of the proposition.\hfill $\Box$

\medskip

Using the bound (\ref{eig}) we can derive from the above proposition that if
the parameters of a strongly regular graph $G$ satisfy $\eta\approx \mu$
then $G$ has a large eigenvalue gap and is therefore a good
pseudo-random graph. We will exhibit several examples of such graphs in
the next section.

\section{Examples}\label{examples}
Here we present some examples of pseudo-random graphs. 
Many of them are well known and already appeared, e.g., in 
\cite{Tho87a} and \cite{Tho87b}, but there also some which have been 
discovered only recently. Since in the rest of the paper we will mostly 
discuss properties of $(n,d,\lambda)$-graphs,
in our examples we emphasize the spectral properties of the constructed 
graphs. We will also use most of these constructions
later to illustrate particular points and to test the strength of 
the theorems.

\noindent
{\bf Random graphs}
\begin{enumerate}
\item
Let $G=G(n,p)$ be a random graph with edge probability $p$. If $p$ satisfies
$pn/\log n \rightarrow \infty$ and $(1-p)n\log n \rightarrow \infty$,
then almost surely all the degrees of $G$ are equal to $(1+o(1))np$.
Moreover it was proved by F\"uredi and Koml\'os \cite{FK}
that the largest eigenvalue of $G$ is a.s. $(1+o(1))np$ and that
$\lambda(G) \leq (2+o(1))\sqrt{p(1-p)n}$. They stated this result only
for constant $p$ but their proof shows that $\lambda(G) \leq O(\sqrt{np})$ 
also when $p\geq poly \log n/n$.
\item
For a positive integer-valued function $d=d(n)$ we define the model 
$G_{n,d}$ of random regular graphs consisting of all regular graphs on 
$n$ vertices of degree $d$ with the uniform probability distribution.
This definition of a random regular graph is conceptually simple, but it 
is not easy to use. Fortunately, for small $d$ there is an efficient
way to generate $G_{n,d}$ which is useful for theoretical studies. This 
is the so called {\em configuration model}. For more details about this 
model, and random regular graphs in general we 
refer the interested reader to two excellent monographs \cite{Bol01}
and \cite{JanLucRuc00}, or to a survey \cite{Wor99}.
As it turns out, sparse random regular graphs have quite different 
properties from those of the binomial random graph $G(n,p), p=d/n$.
For example, they are almost surely 
connected. The spectrum of $G_{n,d}$ for 
a fixed $d$ was studied in \cite{FKS} by 
Friedman, Kahn and Szemer\'edi. Friedman \cite{FRI} proved that 
for constant $d$ the second largest eigenvalue of a random $d$-regular 
graph is $\lambda = (1+o(1))2\sqrt{d-1}$. The approach of 
Kahn and Szemer\'edi gives only $O(\sqrt{d})$ bound on $\lambda$ but 
continues to work also when $d$ is small power of $n$. 
The case $d \gg n^{1/2}$ was recently studied by Krivelevich, Sudakov,
Vu and Wormald \cite{KSVW}. They proved that in this case for any two 
vertices $u,v \in G_{n,d}$ almost surely 
$$\big|codeg(u,v)-d^2/n\big| < Cd^3/n^2 + 6d\sqrt{ \log n} 
/\sqrt{n},$$
where $C$ is some constant and $codeg(u,v)$ is the number of 
common neighbors of $u,v$. Moreover if $d \geq n/\log n$, then $C$ can 
be defined to be zero. Using this it is easy to show that 
for $d \gg n^{1/2}$, the second largest eigenvalue of a random $d$-regular
graph is $o(d)$. The true bound for the second largest 
eigenvalue of $G_{n,d}$ should be probably $(1+o(1))2\sqrt{d-1}$ for 
all values of $d$, but we are still far from proving it.

\noindent
\hspace{-0.95cm}{\bf Strongly regular graphs}
\item
Let $q=p^{\alpha}$ be a prime power which is congruent to $1$ modulo $4$ so
that $-1$ is a square in the finite field $GF(q)$.
Let $P_q$ be the graph whose vertices are all elements of 
$GF(q)$ and two vertices are adjacent if and only if their difference is a 
quadratic residue in $GF(q)$. This graph is usually called 
the {\em Paley graph}. It is easy to see that $P_q$ is $(q-1)/2$-regular.
In addition one can easily compute the number of common neighbors of two 
vertices in $P_q$.
Let $\chi$ be the {\em quadratic residue character} on $GF(q)$, i.e.,
$\chi(0)=0$, $\chi(x)=1$ if $x\not= 0$ and is a square in $GF(q)$ and 
$\chi(x)=-1$ otherwise. By definition, 
$\sum_x\chi(x)=0$ and the number of common neighbors 
of two vertices $a$ and $b$ equals
$$\sum_{x\not=a,b}\left(\frac{1+\chi(a-x)}{2}\right)
\left(\frac{1+\chi(b-x)}{2}\right)=
\frac{q-2}{4}-\frac{\chi(a-b)}{2}+\frac{1}{4}\sum_{x\not=a,b} 
\chi(a-x)\chi(b-x).$$
Using that for $x \not =b$, $\chi(b-x)=\chi\big((b-x)^{-1}\big)$, the last 
term can be rewritten as
$$\sum_{x\not=a,b}\chi(a-x)\chi\big((b-x)^{-1}\big)=
\sum_{x\not=a,b} \chi\Big(\frac{a-x}{b-x}\Big)=
\sum_{x\not=a,b}\chi\Big(1+\frac{a-b}{b-x}\Big)=\sum_{x\not=0,1}\chi(x)=-1.$$
Thus the number of common neighbors of $a$ and $b$ is 
$(q-3)/4-\chi(a-b)/2$. This equals $(q-5)/4$ if $a$ and $b$ are adjacent and 
$(q-1)/4$ otherwise. This implies that 
the Paley graph is a strongly regular graph
with parameters $\big(q, (q-1)/2, (q-5)/4, (q-1)/4\big)$
and therefore its second largest eigenvalue equals
$(\sqrt{q}+1)/2$.

\item
For any odd integer $k$ let $H_k$ denote the graph whose $n_k=2^{k-1}-1$
vertices are all binary vectors of length $k$ with an odd number of
ones except the all one vector,
in which two distinct vertices are adjacent
iff the inner product
of the corresponding vectors is $1$ modulo $2$. Using elementary linear 
algebra it is easy to check that this graph is $(2^{k-2}-2)$-regular.
Also every two nonadjacent vertices vertices in it have 
$2^{k-3}-1$ common neighbors and every two adjacent vertices vertices
have $2^{k-3}-3$ common neighbors. Thus $H_k$ is a strongly regular graph 
with parameters $\big(2^{k-1}-1, 2^{k-2}-2, 2^{k-3}-3, 2^{k-3}-1\big)$ 
and with the second largest eigenvalue 
$\lambda(H_k)=1+2^{\frac{k-3}{2}}$.

\item
Let $q$ be a prime power an let $V(G)$ be the elements of the two 
dimensional vector space over $GF(q)$, so $G$ has $q^2$ vertices.
Partition the $q+1$ lines through the origin of the space into two sets 
$P$ and $N$, where $|P|=k$. Two vertices $x$ and $y$ of the graph $G$
are adjacent if 
$x-y$ is parallel to a line in $P$. This example is due to Delsarte and 
Goethals and to Turyn (see \cite {Sei}). It is easy to check that
$G$ is strongly regular with parameters 
$\big(k(q-1),(k-1)(k-2)+q-2,k(k-1)\big)$. Therefore its eigenvalues, 
besides the trivial one are $-k$ and $q-k$. Thus if $k$ is 
sufficiently large we obtain that $G$ is $d=k(q-1)$-regular graph whose 
second largest eigenvalue is much smaller than $d$.

\noindent
\hspace{-0.95cm}{\bf Graphs arising from finite geometries}
\item
For any integer $t \geq 2$ and for any power $q=p^{\alpha}$ of prime $p$ let
$PG(q,t)$ denote the projective geometry of dimension $t$ over
the finite field $GF(q)$. The interesting case for our purposes here is 
that of large $q$ and fixed $t$. The vertices of $PG(q,t)$ 
correspond to the equivalence classes of the set of all non-zero 
vectors ${\bf x}=(x_0, \ldots, x_t)$ of length $t+1$ over $GF(q)$, where two 
vectors are equivalent if one is a multiple of the other by an element 
of the field. Let $G$ denote the graph whose vertices are the points of 
$PG(q,t)$ and two (not necessarily distinct) vertices ${\bf x}$ and 
${\bf y}$ are adjacent if and only if $x_0y_0+\ldots+x_ty_t=0$. 
This construction is well known. In particular, in case $t=2$ 
this graph is often called the Erd\H os-R\'enyi graph and it contains no 
cycles 
of length $4$. It is easy to see that 
the number of vertices of $G$ is  
$n_{q,t}=\big(q^{t+1}-1\big)/\big(q-1\big)=\big(1+o(1)\big)q^t$
and that it is $d_{q,t}$-regular for 
$d_{q,t}=\big(q^t-1\big)/\big(q-1\big)=\big(1+o(1)\big)q^{t-1}$,
where $o(1)$ tends to zero as $q$ tends to infinity.
It is easy to see that the number of vertices of $G$ with loops is 
bounded by $2\big(q^{t}-1\big)/\big(q-1\big)=\big(2+o(1)\big)q^{t-1}$,
since for every possible value of $x_0, \ldots, x_{t-1}$ we have at most two 
possible choices of $x_t$. Actually using more complicated computation, which 
we omit, one can determine the exact number of vertices with loops.
The eigenvalues of $G$ are easy to compute (see \cite{AK}). Indeed, let 
$A$ be the adjacency matrix of $G$. Then, by the properties of $PG(q,t)$, 
$A^2=AA^T=\mu J+(d_{q,t}-\mu)I$, where 
$\mu=\big(q^{t-1}-1\big)/\big(q-1\big)$, 
$J$ is the all one matrix and $I$ is the identity 
matrix, both of size $n_{q,t} \times n_{q,t}$. Therefore the largest 
eigenvalue of $A$ is $d_{q,t}$ and the absolute 
value of all other eigenvalues is $\sqrt{d_{q,t}-\mu}=q^{(t-1)/2}$. 

\item
The generalized polygons are incidence structures consisting of points 
$\cal P$
and lines $\cal L$. For our purposes we restrict our attention to those 
in which every point is incident to $q+1$ lines and every line is incident to 
$q+1$ points. A generalized $m$-gon defines a bipartite graph $G$ with 
bipartition $({\cal P},{\cal L})$ that 
satisfies the following conditions. The diameter of $G$ is $m$ and for 
every vertex $v \in G$ there is a vertex $u \in G$ such that the 
shortest path from $u$ to $v$ has length $m$. Also for every $r< m$ 
and for every two vertices $u, v$ at distance $r$ there exists a
unique path of length $r$ connecting them. This immediately implies that 
every cycle in $G$ has length at least $2m$. 
For $q \geq 2$, it was proved by Feit and Higman \cite{FH} that
$(q+1)$-regular generalized $m$-gons exist only for $m=3,4,6$.
A {\em polarity} of $G$ is a 
bijection $\pi: {\cal P} \cup {\cal L} \rightarrow {\cal P} \cup {\cal 
L}$ such that $\pi({\cal P})={\cal L}$, $\pi({\cal L})={\cal P}$
and $\pi^2$ is the identity map. Also
for every $p \in {\cal P}, l \in {\cal L}$, $\pi(p)$ is adjacent to 
$\pi(l)$ if and only if $p$ and $l$ are adjacent.
Given $\pi$ we define a polarity 
graph $G^{\pi}$ to be the graph whose vertices are point in $\cal P$ and 
two (not necessarily distinct) points  $p_1, p_2$ are adjacent iff $p_1$ 
was adjacent to $\pi(p_2)$ 
in $G$. Some properties of $G^{\pi}$ can be easily deduced from the 
corresponding properties of $G$. In particular,
$G^{\pi}$ is $(q+1)$-regular and also contains no 
even cycles of length less than $2m$.

For every $q$ which is an odd power of $2$, the incidence graph
of the generalized $4$-gon  has a polarity. The corresponding polarity 
graph
is a $(q+1)$-regular graph with $q^3+q^2+q+1$ vertices. See \cite{BCN},
\cite{LUW} for more details.
This graph contains no cycle of length $6$ and
it is not difficult to compute its eigenvalues (they can be derived,
for example, from
the eigenvalues of the corresponding bipartite incidence graph, given in 
\cite{Ta}).   
Indeed, all the eigenvalues, besides the trivial one (which is $q+1$)
are either $0$ or $\sqrt {2q}$ or $-\sqrt {2q}$.
Similarly, for every $q$ which is an odd power of $3$, the incidence 
graph
of the generalized $6$-gon  has a polarity. The corresponding polarity 
graph
is a $(q+1)$-regular graph with $q^5+q^4+ \cdots +q+1$ vertices (
see again \cite{BCN}, \cite{LUW}).
This graph contains no cycle of length $10$ and
its eigenvalues can be derived using the same technique as in case of 
the $4$-gon. All these eigenvalues, besides the trivial one
are either $\sqrt {3q}$ or $-\sqrt{3q}$ or $\sqrt {q}$ or $-\sqrt {q}$.

\noindent
\hspace{-0.95cm}{\bf Cayley graphs}
\item
Let $G$ be a finite group and let $S$ be a set of non-identity elements 
of $G$ such that $S=S^{-1}$, i.e., for every $s \in S$, $s^{-1}$ also
belongs to $S$. The {\em Cayley graph} $\Gamma(G,S)$ of this group with
respect to  the generating set $S$ is the graph whose set of vertices is $G$ and where 
two vertices  $g$ and $g'$ are adjacent if and only if $g'g^{-1} \in S$. 
Clearly, $\Gamma(G,S)$ is $|S|$-regular and it is connected iff $S$ is a 
set of generators of the group. 
If $G$ is abelian then the eigenvalues of the Cayley graph can be 
computed in terms of the characters of $G$. Indeed, let 
$\chi: G \rightarrow C$ be a character of $G$ and let $A$ be the adjacency 
matrix of $\Gamma(G,S)$ whose rows and columns are indexed by the elements 
of $G$.  Consider the vector ${\bf v}$ defined by ${\bf v}(g)=\chi(g)$. 
Then it is easy to check that $A{\bf v}=\alpha {\bf v}$ with 
$\alpha=\sum_{s\in S}\chi(s)$. In 
addition all eigenvalues can be obtained in this way, since 
every abelian group has exactly $|G|$ different character which are 
orthogonal to each other. Using this fact, one can often give 
estimates on the eigenvalues of $\Gamma(G,S)$ for abelian groups.

One example of a Cayley graph that has already been described earlier 
is $P_q$. In that case the group is the additive group of the finite field 
$GF(q)$ and $S$ is the set of all quadratic residues modulo $q$. Next we 
present a slightly more general construction. Let $q=2kr+1$ be a prime 
power and let $\Gamma$ be a Cayley graph whose group 
is the additive group of $GF(q)$ and whose generating set is
$S=\big\{x=y^k~|~\mbox{for some}~y \in GF(q)\big\}$. 
By definition, $\Gamma$ is $(q-1)/k$-regular. On the other hand,
this graph is not strongly regular unless $k=2$, when it is the Paley  
graph. Let $\chi$ be a nontrivial additive character of $GF(q)$ and 
consider the Gauss sum $\sum_{y \in GF(q)} \chi(y^k)$. Using the 
classical bound $|\sum_{y \in GF(q)} \chi(y^k)|\leq (k-1)q^{1/2}$ 
(see e.g. \cite{LN}) and the above connection between characters and 
eigenvalues we can conclude that the second largest eigenvalue of our 
graph $\Gamma$ is bounded by $O(q^{1/2})$.

\item
Next we present a surprising construction obtained by Alon \cite{A94} of
a  very dense pseudo-random graph that on the other hand is
triangle-free. For a positive integer $k$, consider the 
finite field $GF(2^k)$, whose elements are represented by binary vectors of 
length $k$. If $a, b, c$ are three such vectors, denote by
$(a,b,c)$ the binary vector of length $3k$ whose coordinates are those of 
$a$, followed by coordinates of $b$ and then $c$. Suppose that $k$ is not 
divisible by $3$. Let $W_0$ be the set of all nonzero 
elements $\alpha \in GF(2^k)$ so that the leftmost bit in the binary 
representation of $\alpha^7$ is $0$, and let $W_1$ be the set of all 
nonzero elements $\alpha \in GF(2^k)$ for which the leftmost bit 
of $\alpha^7$ is $1$. Since $3$ does not divide $k$, $7$ does not divide 
$2^k-1$ and hence $|W_0|=2^{k-1}-1$ and $|W_1|=2^{k-1}$, as when
$\alpha$ ranges over all nonzero elements of the field so does 
$\alpha^7$. Let $G_n$ be the graph whose vertices are all
$n=2^{3k}$ binary vectors of length $3k$, where two vectors 
${\bf v}$ and ${\bf v}'$ are adjacent if and only if there exist $w_0\in 
W_0$ and
$w_1\in W_1$ so that ${\bf 
v}-{\bf v}'=(w_0,w_0^3,w_0^5)+(w_1,w_1^3,w_1^5)$, where here 
powers are computed in the field $GF(2^k)$ and the addition is addition 
modulo $2$. Note that $G_n$ is the Cayley graph of the additive group 
${\bf Z}_2^{3k}$ with respect to the generating set 
$S=U_0+U_1$, where $U_0=\big\{(w_0,w_0^3,w_0^5)~|~w_0\in W_0\big\}$
and $U_1$ is defined similarly. A well known fact from Coding 
Theory (see e.g., \cite{MS}), which can be proved using the 
Vandermonde  determinant, is that every set of six distinct vectors in 
$U_0 \cup U_1$ is linearly independent over $GF(2)$. 
In particular all the vectors in $U_0+U_1$ are distinct,
$S=|U_0||U_1|$ and hence $G_n$ is $|S|=2^{k-1}(2^{k-1}-1)$-regular.
The statement that $G_n$ is triangle free is clearly equivalent to the 
fact
that the sum modulo $2$ of any set of $3$ nonzero elements of 
$S$ is not a zero-vector. Let $u_0+u_1, u'_0+u'_1$ and $u''_0+u''_1$
be three distinct element of $S$, where $u_0,u'_0,u''_0 \in U_0$ 
and $u_1,u'_1,u''_1 \in U_1$. By the above discussion, if the sum of these 
six vectors is zero, then every vector must appear an even number of times
in the sequence $(u_0,u'_0,u''_0,u_1,u'_1,u''_1)$. However, since $U_0$ and 
$U_1$ are disjoint, this is clearly impossible. Finally, as we already 
mentioned, the eigenvalues of $G_n$ can be computed in terms of characters 
of 
${\bf Z}_2^{3k}$. Using this fact together with the Carlitz-Uchiyama bound 
on the characters of ${\bf Z}_2^{3k}$ it was proved in \cite{A94} that the 
second eigenvalue of $G_n$ is bounded by $\lambda \leq 9\cdot 
2^k+3\cdot 2^{k/2}+1/4$.  

\item
The construction above can be extended in the obvious way as mentioned in 
\cite{ALONK}. 
Let $h\geq 1$ and suppose that $k$ is an integer such that 
$2^k-1$ is not divisible by $4h+3$.
Let $W_0$ be the set of all nonzero
elements $\alpha \in GF(2^k)$ so that the leftmost bit in the binary
representation of $\alpha^{4h+3}$ is $0$, and let $W_1$ be the set of all
nonzero elements $\alpha \in GF(2^k)$ for which the leftmost bit
of $\alpha^{4h+3}$ is $1$. Since $4h+3$ does not divide
$2^k-1$ we have that $|W_0|=2^{k-1}-1$ and $|W_1|=2^{k-1}$, as when
$\alpha$ ranges over all nonzero elements of the field so does
$\alpha^{4h+3}$.
Define $G$ to be the Cayley graph of the additive group
${\bf Z}_2^{(2h+1)k}$ with respect to the generating set
$S=U_0+U_1$, where $U_0=\big\{(w_0,w_0^3,\ldots,w_0^{4h+1})~|~w_0\in 
W_0\big\}$ and $U_1$ is defined similarly. Clearly,
$G$ is a $2^{k-1}(2^{k-1}-1)$-regular graph on $2^{(2h+1)k}$ vertices. 
Using methods from \cite{A94}, one can show that
$G$ contains no odd cycle of length $\leq 2h+1$ and that the 
second eigenvalue of $G$ is bounded by $O(2^k)$.

\item
Now we describe the celebrated expander graphs constructed by 
Lubotzky, Phillips and Sarnak \cite{LPS} and independently by Margulis 
\cite{Margulis}.
Let $p$ and $q$ be unequal primes, both congruent to $1$ modulo $4$
and such that $p$ is a quadratic residue modulo $q$. As usual denote by
$PSL(2,q)$ the factor group of the group of two by two matrices over
$GF(q)$ with determinant $1$ modulo its normal subgroup consisting of the 
two 
scalar matrices $\bigg(\begin{array}{cc}1&0\\0&1\end{array}\bigg)$ and
$\bigg(\begin{array}{cc}-1&0\\0&-1\end{array}\bigg)$. The 
graphs we describe are Cayley graphs of 
$PSL(2,q)$. A well known theorem of Jacobi asserts that the number of 
ways to represent a positive integer $n$ as a sum of $4$ squares is 
$8\sum_{4 \not \, |\, d, d|n}d$. This easily implies that there are precisely 
$p+1$ 
vectors ${\bf a}=(a_0, a_1, a_2, a_3)$, where $a_0$ is an odd positive 
integer, $a_1, a_2, a_3$ are even integers  and $a_0^2+ a_1^2+a_2^2+a_3^2=p$.
From each such vector construct the matrix $M_a$ in $PSL(2,q)$ where
$M_a=\frac{1}{\sqrt{p}}
\bigg(\begin{array}{cc}a_0+ia_1&a_2+ia_3\\-a_2+ia_3&a_0-ia_1\end{array}\bigg)$
and $i$ is an integer satisfying $i^2=-1(\mbox{mod}~ q)$.
Note that, indeed, the determinant of $M_a$ is $1$ and that the square root 
of $p$ modulo $q$ does exist.
Let $G^{p,q}$ denote the Cayley graph of $PSL(2,q)$ with respect to these
$p+1$ matrices. In \cite{LPS} it was proved that
if $q >2\sqrt{p}$ then $G^{p,q}$ is a connected $(p+1)$-regular graph on 
$n=q(q^2-1)/2$ vertices. Its girth is at least $2\log_p q$ and all the 
eigenvalues of its adjacency matrix, besides the trivial one
 $\lambda_1=p+1$, are at most $2 \sqrt{p}$ in absolute value.
 The bound on the eigenvalues was obtained by applying 
deep results of Eichler and Igusa concerning the Ramanujan conjecture.
The graphs $G^{p,q}$ have very good expansion properties and have numerous 
applications in Combinatorics and Theoretical Computer Science.

\item
The {\em projective norm graphs} $NG_{p,t}$ have been constructed in
\cite{ARS}, modifying an earlier construction given in \cite{KRS}.
These graphs {\bf are not} Cayley graphs, but as one will immediately see, 
their construction has a similar flavor.
The construction is the following. Let $t>2$ be an integer, let $p$
be a prime, let $GF(p)^*$ be  the multiplicative group of the
field with $p$ elements
and let $GF(p^{t-1})$ be the field with  $p^{t-1}$ elements.
The set of vertices of the graph $NG_{p,t}$ is the set
$V=GF(p^{t-1})\times GF(p)^*$. Two distinct vertices $(X,a)$
and $(Y,b)\in V$ are adjacent if and only if $N(X+Y)=ab$, where the norm
$N$ is understood over $GF(p)$, that is, $N(X)=X^{1+p+\cdots+p^{t-2}}.$
Note that $|V|=p^t-p^{t-1}$. If $(X,a)$ and $(Y,b)$ are adjacent,
then $(X,a)$ and $Y\neq -X$ determine $b$.
Thus $NG_{p,t}$ is a regular graph of degree $p^{t-1}-1$. In addition, it was 
proved in \cite{ARS}, that $NG_{p,t}$ contains no complete bipartite graphs 
$K_{t,(t-1)!+1}$.
These graphs can be also defined in the same manner starting with a 
prime power instead of the prime $p$. It is also not difficult to 
compute the eigenvalues of this graph. Indeed, put $q=p^{t-1}$ and let 
$A$ be the adjacency matrix of $NG_{p,t}$.
The rows and columns of this matrix are indexed by the ordered pairs
of the set $GF(q) \times GF(p)^*$. Let $\psi$ be a character of the
additive group of $GF(q)$, and let $\chi$ be a character of the
multiplicative group of $GF(p)$. Consider the vector
${\bf v}: GF(q) \times GF(p)^* \mapsto C$ defined by
${\bf v}(X,a)=\psi(X) \chi(a)$. Now one can check (see \cite{AR}, 
\cite{Sz} for more 
details) that the vector 
${\bf v}$ is an eigenvector of $A^2$ with eigenvalue 
$\big| \sum_{Z\in GF(q), Z\not=0} \psi(Z)\chi(N(Z))\big|^2$ and that all
eigenvalues of $A^2$ have this form. 
Set $\chi'(Z)=\chi(N(Z))$ for all nonzero $Z$ in $GF(q)$. Note that as the
norm is multiplicative, $\chi'$ is a multiplicative
character of the large field. Hence the above expression
is a square of the absolute value of the Gauss sum and it is well known (see 
e.g. \cite{Da}, 
\cite{Bol01}) that the value of each such square, besides the trivial one 
(that is,
when either $\psi$ or $\chi'$ are trivial),  is $q$.
This implies that the second largest eigenvalue of $NG_{p,t}$
is $\sqrt{q}=p^{(t-1)/2}$.
\end{enumerate}

\section{Properties of pseudo-random graphs}
We now examine closely properties of pseudo-random graphs, with a
special emphasis on $(n,d,\lambda)$-graphs. The majority of them
are obtained using the estimate (\ref{eig}) of Theorem \ref{eigen}, 
showing
again the extreme importance and applicability of the latter result. It
is instructive to compare the properties of pseudo-random graphs,
considered below, with the analogous properties of random graphs, usually
shown to hold by completely different methods. The set of properties we
chose to treat here is not meant to be comprehensive or systematic, but
quite a few rather diverse graph parameters will be covered.

\subsection{Connectivity and perfect matchings}
The {\em vertex-connectivity} of a graph $G$ is the minimum number of 
vertices that 
we need to delete to make $G$ disconnected. We denote this parameter by
$\kappa(G)$. For random graphs it is well known (see, e.g., \cite{Bol01}) 
that the vertex-connectivity is almost surely the same as the minimum 
degree.
Recently it was also proved (see \cite{KSVW} and \cite{CFR}) that 
random $d$-regular graphs are $d$-vertex-connected. For 
$(n,d,\lambda)$-graphs it 
is easy to show the following.

\begin{theo}
\label{connectivity}
Let $G$ be an $(n,d,\lambda)$-graph with $d \leq n/2$.
Then the vertex-connectivity of $G$ satisfies:
$$\kappa(G) \geq d-36\lambda^2/d.$$
\end{theo}

\noindent
{\bf Proof.}\, We can assume that $\lambda \leq d/6$, since otherwise 
there is nothing to prove. Suppose that there is a subset $S \subset V$ of 
size less than $d-36\lambda^2/d$
such that the induced graph $G[V-S]$ is disconnected. Denote by $U$ the
set of vertices of the smallest connected component of $G[V-S]$ and set
$W=V-(S \cup U)$. Then $|W| \geq (n-d)/2\geq n/4$ and there is no edge 
between $U$ and $W$. Also $|U|+|S|>d$, since all the neighbors of 
a vertex 
from $U$ are contained in $S \cup U$. Therefore $|U|\geq 36\lambda^2/d$. 
Since there are no edges between $U$ and $W$, by Theorem \ref{eigen},
we have that $d|U||W|/n <\lambda \sqrt{|U||W|}$. This implies that
$$|U|<\frac{\lambda^2n^2}{d^2|W|}= \frac{\lambda}{d}\frac{n}{|W|} 
\frac{\lambda n}{d} \leq \frac{1}{6} \cdot 4 \cdot\frac{\lambda n}{d}<
\frac{\lambda n}{d}.$$
Next note that, by Theorem \ref{eigen}, the number of edges spanned by 
$U$ is at most
$$e(U) \leq \frac{d|U|^2}{2n}+\frac{\lambda|U|}{2}<
\frac{\lambda n}{d}\frac{d|U|}{2n}+\frac{\lambda|U|}{2}=
\frac{\lambda|U|}{2}+\frac{\lambda|U|}{2}= \lambda|U|.$$
As the degree of every vertex in $U$ is $d$, it follows that
$$e(U,S)\geq d|U|-2e(U)> (d-2\lambda)|U|\geq 2d|U|/3.$$
On the other hand 
using again Theorem \ref{eigen} together with the facts that 
$|U|\geq 36\lambda^2/d$, $|S|<d$ and $d \leq n/2$ we conclude that
\begin{eqnarray*}
e(U,S)&\leq& \frac{d|U||S|}{n}+\lambda\sqrt{|U||S|}
<\frac{d}{n}d|U|+\lambda\sqrt{d|U|} \leq \frac{d|U|}{2}+
\frac{\lambda\sqrt{d}|U|}{\sqrt{|U|}}\\
&\leq&\frac{d|U|}{2}+\frac{\lambda\sqrt{d}|U|}{6\lambda/\sqrt{d}}=
\frac{d|U|}{2}+\frac{d|U|}{6}=\frac{2d|U|}{3}.
\end{eqnarray*}
This contradiction completes the proof.\hfill $\Box$

\medskip

\noindent
The constants in this theorem can be easily improved and we make no 
attempt to optimize them. Note that, in particular, for an
$(n,d,\lambda)$-graph 
$G$ with $\lambda=O(\sqrt{d})$ we have that $\kappa(G)= d-\Theta(1)$.

Next we present an example which shows that the assertion of Theorem
\ref{connectivity} is tight up to a constant factor.
Let $G$ be any $(n,d,\lambda)$-graph with $\lambda=\Theta(\sqrt{d})$.
We already constructed several such graphs in the previous section. For an 
integer $k$, consider a new graph $G_k$, which is obtained by 
replacing each vertex of $G$ by the complete graph of order $k$ and by
connecting two vertices of $G_k$ by an edge if and only if the 
corresponding vertices of $G$ are connected by an edge.
Then it follows immediately from the definition that 
$G_k$ has $n'=nk$ vertices  and is $d'$-regular graph with $d'=dk+k-1$.
Let $\lambda'$ be the second eigenvalue of $G_k$. To estimate 
$\lambda'$ note that the adjacency matrix of $G_k$ equals to
$A_G\otimes J_k+I_n\otimes A_{K_k}$. Here $A_G$ is the adjacency matrix of 
$G$, $J_k$ is the all one matrix of size $k\times k$, $I_n$ is the 
identity
matrix of size $n \times n$ and $A_{K_k}$ is the adjacency matrix of the 
complete graph of order $k$. Also the tensor product of the $m\times n$ 
dimensional matrix $A=(a_{ij})$ and the $s\times t$-dimensional matrix 
$B=(b_{kl})$ is the $ms\times nt$-dimensional matrix $A\otimes B$, whose 
entry labeled $((i,k)(j,l))$ is $a_{ij}b_{kl}$. In case $A$ and $B$ are 
symmetric matrices with spectrums $\{\lambda_1, \ldots , \lambda_n\}$, $\{ 
\mu_1,\ldots , \mu_t\}$ respectively, it is a simple consequence of the 
definition that the spectrum of $A\otimes B$ is $\{ \lambda_i\mu_k: i=1, 
\ldots ,n, k=1,\ldots , t\}$ (see, e.g. \cite{Lov93}). 
Therefore the second eigenvalue of $A_G\otimes J_k$ is $k\lambda$.
On the other hand $I_n\otimes A_{K_k}$ is the adjacency matrix of the 
disjoint union of $k$-cliques and therefore the absolute value of all its 
eigenvalues is at most $k-1$.
Using these two facts we conclude that $\lambda'\leq \lambda k+k-1$ and 
that $G_k$ is $(n'=nk,d'=dk+k-1,\lambda'=\lambda k+k-1)$-graph. 
Also it is easy to see that the set of vertices 
of $G_k$ that corresponds to a vertex in $G$ has exactly $dk$  
neighbors outside this set. By deleting these neighbors we can disconnect 
the graph $G_k$ and thus 
$$\kappa(G_k) \leq dk=d'-(k-1)=d'-\Omega\big((\lambda')^2/d'\big).$$

Sometimes we can improve the result of Theorem \ref{connectivity}
using the information about co-degrees of vertices in our graph.  
Such result was used in \cite{KSVW} to determine the vertex-connectivity 
of dense random $d$-regular graphs.

\begin{prop}
\label{p42}\cite{KSVW}
Let $G=(V,E)$ be a
$d$-regular graph on $n$ vertices such that $ \sqrt{n}\log n < d \leq 
3n/4$ and the number of common neighbors for every two 
distinct vertices in $G$ is $(1+o(1))d^2/n$. Then the graph $G$ is
$d$-vertex-connected. 
\end{prop}

Similarly to vertex-connectivity, define the {\em edge-connectivity} of 
a graph $G$ to be the minimum number of
edges that we need to delete to make $G$ disconnected. We denote this
parameter by $\kappa'(G)$. Clearly the edge-connectivity is always at most
the minimum degree 
of a graph. We also say that $G$ has a {\em perfect matching} if there
is a set of disjoint edges that covers all the vertices of $G$.
Next we show that $(n,d,\lambda)$-graphs even with a very weak 
spectral gap are $d$-edge-connected and have a perfect matching 
(if the number of vertices is even).

\begin{theo}
\label{edge-connectivity}
Let $G$ be an $(n,d,\lambda)$-graph with $d-\lambda\geq 2$.
Then $G$ is $d$-edge-connected. When $n$ is even,  
it has a perfect matching.
\end{theo}

\noindent
{\bf Proof.}\, 
Let $U$ be a subset of vertices of $G$ of size at most $n/2$. To prove 
that $G$ is $d$-edge-connected we need to show that 
there are always at least $d$ edges between $U$ and $V(G)-U$.
If $1 \leq |U|\leq d$, then every vertex in $U$ has at least $d-(|U|-1)$ 
neighbors outside $U$ and therefore $e(U,V(G)-U) \geq |U|\big(d-|U|+1\big) 
\geq d$. On the other hand if $d \leq |U|\leq n/2$, then using that
$d-\lambda\geq 2$ together with Theorem 
\ref{eigen} we obtain that
\begin{eqnarray*}
e\big(U,V(G)-U\big) &\geq& 
\frac{d|U|(n-|U|)}{n}-\lambda\sqrt{|U|(n-|U|)\left(1-\frac{|U|}{n}\right)
\left(1-\frac{n-|U|}{n}\right)}\\
&=&(d-\lambda)\frac{(n-|U|)}{n}|U| \geq
2\cdot\frac{1}{2}\cdot|U|=|U|\geq d,
\end{eqnarray*}
and therefore $\kappa'(G)=d$. 

To show that $G$ contains a perfect matching we apply the celebrated 
Tutte's condition. Since $n$ is even, we need to prove that for every 
nonempty set of vertices 
$S$, the 
induced graph $G[V-S]$ has at most $|S|$ connected components of odd size.
Since $G$ is $d$-edge-connected we have that there are at least $d$ edges 
from every connected component of 
$G[V-S]$ to $S$. On the other hand there are at most $d|S|$ edges incident 
with vertices in $S$. Therefore $G[V-S]$ has at most $|S|$ connected 
components and hence $G$ contains a perfect matching. 
\hfill $\Box$

\subsection{Maximum cut}
Let $G=(V,E)$ be a graph and let $S$ be a nonempty proper
subset of $V$. Denote by $(S,V-S)$ the cut of $G$ consisting of all
edges with one end in $S$ and another one in $V-S$. The {\em size}
of the cut is the number of edges in it. The MAX CUT problem is the
problem of finding a cut of maximum size in $G$. Let $f(G)$ be the 
size of the  maximum cut in $G$. 
MAX CUT is one of the most natural combinatorial
optimization problems. It is well known that this problem is 
NP-hard \cite{GJ}. Therefore  it is useful to have bounds on 
$f(G)$ based on other parameters of the graph, that can be computed 
efficiently. 

Here we describe two such folklore results. First, consider a random partition 
$V=V_1\cup V_2$, obtained by assigning each vertex $v\in V$ to $V_1$ or $V_2$ 
with probability $1/2$ independently.
It is easy to see that each edge of $G$ has probability $1/2$ to cross
between $V_1$ and $V_2$. Therefore the expected
number of edges in the cut $(V_1,V_2)$ is $m/2$, where $m$ is the number of edges 
in $G$. This implies that for every graph $f(G)\geq m/2$. The example of
a complete graph shows that this lower bound is asymptotically optimal.
The second result provides an upper bound for $f(G)$, for a regular
graph $G$, in terms of the smallest eigenvalue of its adjacency matrix.

\begin{prop}
\label{max-cut}
Let $G$ be a $d$-regular graph (which may have loops) of order $n$ with
$m=dn/2$ edges and let
$\lambda_1\geq \lambda_2 \geq \ldots \geq \lambda_n$ be the eigenvalues of the
adjacency matrix of $G$. Then
$$f(G) \leq \frac{m}{2}-\frac{\lambda_n n}{4}.$$
In particular if $G$ is an $(n,d,\lambda)$-graph then 
$f(G) \leq (d+\lambda)n/4$. 
\end{prop}   

\noindent
{\bf Proof.}\, Let $A=(a_{ij})$ be the adjacency matrix of
$G=(V,E)$ and let $V=\{1, \ldots, n\}$. Let
${\bf x}=(x_1, \ldots, x_n)$ be any vector with coordinates $\pm 1$.
Since the graph $G$ is $d$-regular we have
$$\sum_{(i,j)\in E} (x_i-x_j)^2=d\sum_{i=1}^n x_i^2-
\sum_{i,j}a_{ij}x_ix_j=dn-{\bf x}^tA{\bf x}.$$
 By the variational definition of the eigenvalues of $A$,
for any vector $z \in R^n$, $z^tAz \geq \lambda_n\|z\|^2$.
Therefore
\begin{equation}
\label{a}
\sum_{(i,j)\in E} (x_i-x_j)^2=dn-{\bf x}^tA{\bf x} \leq
dn-\lambda_n\|{\bf x}\|^2=dn-\lambda_nn.
\end{equation}

Let $V=V_1\cup V_2$ be an arbitrary partition of $V$ into two disjoint
subsets and let $e(V_1,V_2)$ be the number of edges in the bipartite
subgraph of $G$ with bipartition $(V_1,V_2)$. For every vertex $v \in
V(G)$ define $x_v=1$ if $v \in V_1$ and $x_v=-1$ if $v \in V_2$.
Note that for every edge $(i,j)$ of $G$,
$(x_i-x_j)^2=4$ if this edge has its ends
in the distinct parts of the above partition and is zero otherwise.
Now using (\ref{a}), we conclude that
$$\hspace{3.2cm}
e(V_1,V_2)=\frac{1}{4} \sum_{(i,j)\in E} (x_i-x_j)^2 \leq
\frac{1}{4}(dn-\lambda_nn)=\frac{m}{2}-\frac{\lambda_nn}{4}.
\hspace{3.2cm}\Box$$

This upper bound is often used to show that some particular results
about maximum cuts are tight. For example this approach was used in
\cite{ALON} and \cite{ABKS}. In these papers the authors proved that
 for every graph $G$ with $m$ edges and girth at least $r\geq 4$,
$f(G) \geq m/2 +\Omega\big(m^{\frac{r}{r+1}}\big)$.
They also show, using Proposition \ref{max-cut} and Examples 9, 6 from Section 3, 
that this bound is tight for $r=4,5$.

\subsection{Independent sets and the chromatic number}
The {\em independence number} $\alpha(G)$ of a graph $G$
is the maximum cardinality of a set of vertices of $G$ no two of which are 
adjacent. Using Theorem \ref{eigen} we can immediately establish an upper 
bound on the size of a maximum independent set of pseudo-random graphs.
\begin{prop}
\label{ind-set}
Let $G$ be an $(n,d,\lambda)$-graph, then
$$\alpha(G)\leq \frac{\lambda n}{d+\lambda}.$$
\end{prop}

\noindent
{\bf Proof.}\, Let $U$ be an independent set in $G$, then 
$e(U)=0$ and by Theorem  \ref{eigen} we have that
$d|U|^2/n\leq \lambda |U|(1-|U|/n)$. This implies that
$|U|\leq \lambda n/(d+\lambda)$. \hfill $\Box$

\vspace{0.15cm}
\noindent
Note that even when $\lambda=O(\sqrt{d})$ this bound only has order of magnitude
$O(n/\sqrt{d})$. This contrasts sharply
with the behavior of random graphs where it is known 
(see \cite{Bol01} and \cite{JanLucRuc00}) that the independence number of random graph 
$G(n,p)$ is only $\Theta\big(\frac{n}{d}\log d\big)$ where
$d=(1+o(1))np$. More strikingly there are graphs for which 
the bound in Proposition \ref{ind-set} cannot be improved. One such
graph is the 
Paley graph $P_q$ with
$q=p^2$ (Example 3 in the previous section). 
Indeed it is easy to see that 
in this case all elements of the subfield $GF(p)\subset GF(p^2)$ are quadratic residues in 
$GF(p^2)$. This implies that for every quadratic non-residue $\beta \in GF(p^2)$ 
all elements of any multiplicative coset $\beta GF(p)$ form an independent 
set of size 
$p$. As we already mentioned, $P_q$ is an 
$(n,d,\lambda)$-graph with $n=p^2,d=(p^2-1)/2$ and $\lambda=(p+1)/2$. 
Hence for this graph we get $\alpha(P_q)=\lambda n/(d+\lambda)$.

Next we obtain a lower bound on the independence number 
of pseudo-random graphs.  We present  a slightly 
more general result by Alon et al. \cite{AKS} which we will need later.
\begin{prop}
\label{ind-set1}\cite{AKS}
Let $G$ be an $(n,d,\lambda)$-graph such that 
$\lambda<d\leq  0.9n$. 
Then the induced subgraph $G[U]$ of $G$ on any
subset $U, |U|=m$, contains an independent set of size at least
$$\alpha(G[U]) \geq \frac{n}{2(d-\lambda)}
\ln \left(\frac{m(d-\lambda)}{n(\lambda+1)} +1 \right).$$
In particular,
$$\alpha(G) \geq \frac{n}{2(d-\lambda)}
\ln \left(\frac{(d-\lambda)}{(\lambda+1)} +1 \right).$$
\end{prop}

\noindent
{\bf Sketch of proof.}\,
First using Theorem \ref{eigen} it is easy to show that if
$U$ is a set of $bn$ vertices of $G$, then the minimum degree in the
induced subgraph $G[U]$ is at most $db+\lambda(1-b)=(d-\lambda)b+\lambda$.
Construct an independent set $I$ in the induced subgraph
$G[U]$ of $G$ by the following greedy procedure.
Repeatedly choose a vertex of minimum degree
in $G[U]$  , add it to the independent set $I$ and delete it and its
neighbors from $U$, stopping
when the remaining set of vertices is empty. 
Let $a_i, i\geq 0$ be the sequence of numbers defined by the
following recurrence formula:  
$$a_0=m,~~ 
a_{i+1}=a_i-\left(d\frac{a_i}{n}+\lambda(1-\frac{a_i}{n})+1\right)=
\left(1-\frac{d-\lambda}{n}\right)a_i-(\lambda+1),~\forall i \geq 0.$$
By the above discussion, it is easy to see
that the size of the remaining set of vertices after $i$ iterations is at 
least $a_i$. Therefore
the size of the resulting independent set $I$ is at least the smallest
index $i$ such that $a_i\leq 0$. By solving the recurrence equation we 
obtain that this index satisfies:
$$\hspace{5cm}i \geq \frac{n}{2(d-\lambda)}
\ln \left(\frac{m(d-\lambda)}{n(\lambda+1)} +1 \right)\,. \hspace{5cm} 
\Box$$

\vspace{0.15cm}
\noindent
For an $(n,d,\lambda)$-graph $G$ with $\lambda \leq d^{1-\delta},
\delta>0$, this proposition 
implies that $\alpha(G) \geq \Omega\big(\frac{n}{d}\log d\big)$. This shows that
the independence number of a
pseudo-random graph with a sufficiently small second eigenvalue is up to a 
constant factor 
at least as large as $\alpha(G(n,p))$ with $p=d/n$. On the other hand the graph $H_k$
(Example 4, Section 3) 
shows that even when $\lambda \leq O(\sqrt{d})$ the independence number
of $(n,d,\lambda)$-graph can be smaller than $\alpha(G(n,p))$ with 
$p=d/n$.
This graph has $n=2^{k-1}-1$ vertices, degree $d=(1+o(1))n/2$ and $\lambda=
\Theta(\sqrt{d})$. Also it is easy to see that every independent set in 
$H_k$ corresponds to a family of orthogonal vectors in ${\bf Z}_2^k$ and thus has size at 
most $k=(1+o(1))\log_2 n$. This is only half of the size of a maximum
independent set in the corresponding random graph $G(n,1/2)$.

A {\em vertex-coloring} of a graph $G$ is an assignment  of a color
to each of its vertices. The coloring is {\em proper} if no two
adjacent vertices get the same color. The {\em chromatic number}
$\chi(G)$ of $G$ is the minimum number of colors used in a proper coloring of
it.  Since every color class in the proper coloring of $G$ forms an independent 
set we can immediately obtain that $\chi(G)\geq |V(G)|/\alpha(G)$. This together 
with Proposition \ref{ind-set} implies the following result of Hoffman \cite{Ho}.

\begin{coro}
\label{hofman}
Let $G$ be an $(n,d,\lambda)$-graph. Then the chromatic number of $G$ is at least
$1+d/\lambda$.
\end{coro}

On the other hand, using Proposition \ref{ind-set1}, one can obtain the following 
upper bound on the chromatic number of pseudo-random graphs.

\begin{theo}
\label{chromatic}\cite{AKS}
Let $G$ be an $(n,d,\lambda)$-graph
such that $\lambda<d\leq 0.9n$.
Then the chromatic number of $G$ satisfies
$$
\chi(G) \leq \frac{6(d-\lambda)}{\ln\big (\frac{d-\lambda}{\lambda+1} 
+1\big)}\, .
$$
\end{theo}

\noindent
{\bf Sketch of proof.}\,
Color the graph $G$ as follows. As long as the remaining 
set of vertices $U$ contains at least
$n/\ln (\frac{d-\lambda}{\lambda+1}+1)$ vertices, by Proposition
\ref{ind-set1} we can find an 
independent
set of vertices in the induced subgraph $G[U]$ of size at least
$$\frac{n}{2(d-\lambda)} \ln \left(
\frac{|U|(d-\lambda)}{n(\lambda+1)}+1\right) \geq
\frac{n}{4(d-\lambda)} \ln \left(
\frac{d-\lambda}{\lambda+1}+1\right).$$
Color all the members of such a set
by a new color, delete them from the graph and   
continue.
When this process terminates, the remaining set of vertices $U$ is
of size at most $n/\ln (\frac{d-\lambda}{\lambda+1}+1)$ and we used at
most $4(d-\lambda)/\ln(\frac{d-\lambda}{\lambda+1}+1)$ colors so far.
As we already mentioned above, for every subset $U'\subset U$ the induced
subgraph $G[U']$ contains a vertex of degree at most
$$(d-\lambda)\frac{|U'|}{n}+\lambda\leq
(d-\lambda)\frac{|U|}{n}+\lambda\leq
\frac{d-\lambda}{\ln (\frac{d-\lambda}{\lambda+1}+1)}+\lambda \leq
\frac{2(d-\lambda)}{\ln (\frac{d-\lambda}{\lambda+1}+1)}-1.$$
Thus we
can complete the coloring of $G$ by
coloring $G[U]$ using
at most $2(d-\lambda)/\ln (\frac{d-\lambda}{\lambda+1}+1)$
additional colors.
The total number of colors used is at most
$6(d-\lambda)/\ln (\frac{d-\lambda}{\lambda+1}+1)$. \hfill$\Box$

\vspace{0.15cm}
\noindent
For an $(n,d,\lambda)$-graph $G$ with $\lambda \leq d^{1-\delta}, \delta>0$ this proposition
implies that $\chi(G) \leq O\big(\frac{d}{\log d}\big)$. This shows that
the chromatic number of a
pseudo-random graph with a sufficiently small second eigenvalue is up to a 
constant factor
at least as small as $\chi(G(n,p))$ with $p=d/n$. 
On the other hand, the 
Paley graph $P_q, q=p^2$, shows that sometimes the chromatic number of 
a pseudo-random graph can be much smaller than the above bound, even in 
the case 
$\lambda=\Theta(\sqrt{d})$. Indeed, as we already mentioned above, 
all elements of the subfield $GF(p)\subset GF(p^2)$ are quadratic residues in
$GF(p^2)$. This implies that for every quadratic non-residue $\beta \in GF(p^2)$
all elements of a multiplicative coset $\beta GF(p)$ form an independent 
set of size
$p$. Also all additive cosets of $\beta GF(p)$ are independent sets in
$P_q$. This implies that $\chi(P_q)\leq \sqrt{q}=p$. In fact $P_q$ contains a clique of size 
$p$ (all elements of a subfield $GF(p)$), showing that
$\chi(P_q)=\sqrt{q}\ll q/\log q$. Therefore the bound in Corollary 
\ref{hofman} is best possible.

A more complicated quantity related to the chromatic number is the  {\em 
list-chromatic 
number}
$\chi_l(G)$ of $G$, introduced in \cite{ERT} and \cite{Vi}.  This is the
minimum integer $k$ such that for every assignment of a set $S(v)$ of
$k$ colors to every vertex $v$ of $G$, there is a proper coloring of  
$G$ that assigns to each vertex $v$ a color from $S(v)$. The study of this 
parameter received a considerable amount of attention in
recent years, see, e.g., \cite{Al2}, \cite{KTV} for two surveys.
Note that from the definition it follows immediately that $\chi_l(G) \geq 
\chi(G)$ and it is known that the gap between these two parameters
can be arbitrarily large. The list-chromatic number of pseudo-random 
graphs was 
studied by Alon, Krivelevich and Sudakov \cite{AKS} and independently by Vu 
\cite{Vu}. In \cite{AKS} and \cite{Vu} the authors mainly considered graphs with all degrees 
$(1+o(1))np$ and all co-degrees $(1+o(1))np^2$. Here we use ideas from these two 
papers to obtain an upper bound on the list-chromatic number of 
an $(n,d,\lambda)$-graphs.
This bound has the same order of magnitude as the list chromatic number of
the truly random 
graph $G(n,p)$ with $p=d/n$ (for more details see \cite{AKS}, \cite{Vu}).

\begin{theo}
\label{choice}
Suppose that $0<\delta<1$ and let $G$ be an $(n,d,\lambda)$-graph satisfying
$\lambda \leq d^{1-\delta}$, $d\le 0.9n$. 
Then the list-chromatic number of $G$ is bounded by
$$\chi_l(G)\leq O\left(\frac{d}{\delta\log d}\right).$$ 
\end{theo}

\noindent
{\bf Proof.}\, Suppose that $d$ is sufficiently large and consider first the 
case when $d \leq n^{1-\delta/4}$. Then by 
Theorem \ref{eigen} the neighbors of every vertex in $G$ span
at most $d^3/n+\lambda d \leq O(d^{2-\delta/4})$ edges. Now we can apply the 
result of Vu \cite{Vu} which says that if the neighbors of every vertex in 
a graph $G$ with maximum  degree $d$ span at most $O(d^{2-\delta/4})$ 
edges then 
$\chi_l(G)\leq O\big(d/(\delta\log d)\big).$ 

Now consider the case when $d \geq n^{1-\delta/4}$. For every 
vertex $v \in V$, let $S(v)$ be a list of at least $\frac{7d}{\delta \log n}$ 
colors.
Our objective is to prove that there is a proper coloring of $G$
assigning to each vertex a color from its list. As long as there is
a set $C$ of at least
$n^{1-\delta/2} $ vertices containing the same color  $c$ in their
lists we can, by
Proposition \ref{ind-set1}, find an independent set of at least
$\frac{\delta n}{6d} \log n$ vertices in $C$, color them all by
$c$,  omit them from the graph and  omit the color $c$ from all
lists. The total number of colors that can be deleted in this
process cannot exceed $\frac{6d}{\delta \log n}$ (since in each
such deletion at least $\frac{\delta n}{6d}\log n$ vertices are
deleted from the graph). When this process terminates, no color
appears in more than $n^{1-\delta/2}$ lists, and each list still
contains at least $\frac{d}{\delta \log n} > n^{1-\delta/2}$     
colors. Therefore, by Hall's theorem, we can assign to each of the
remaining vertices
a color from its list so that no color is being assigned to more
than one vertex,
thus completing the coloring and the proof.\hfill $\Box$

\subsection{Small subgraphs}
We now examine small subgraphs of pseudo-random graphs. 
Let $H$ be a fixed graph of order $s$ with $r$ edges and with automorphism group
$Aut(H)$. Using the 
second moment method it is 
not difficult to show that for every constant $p$ the random graph $G(n,p)$ contains 
$$(1+o(1))p^r(1-p)^{{s\choose 2}-r}\frac{n^s}{|Aut(H)|}$$
induced copies of $H$. Thomason extended 
this result to jumbled graphs. He showed in \cite{Tho87a} that if a graph 
$G$ is 
$(p,\alpha)$-jumbled and $p^sn\gg 42\alpha s^2$ then the number of induced subgraphs of $G$ 
which are isomorphic to $H$ is
$(1+o(1))p^s(1-p)^{{s\choose 2}-r}n^s/|Aut(H)|$.

Here we present a result of Noga Alon \cite{Alon} that
proves that every large subset of the set of vertices
of $(n,d,\lambda)$-graph contains the "correct" number of copies of any
fixed sparse graph. An additional advantage of this result is that its assertion depends
not on the number of vertices $s$ in $H$ but only on its maximum degree
$\Delta$ which can be smaller than $s$.
Special cases of this result have appeared in various
papers including \cite{AK}, \cite{AP} and probably other papers as well.
The approach here is similar to the one in \cite{AP}.
\begin{theo}
\label{number-subgraphs} \cite{Alon}
Let $H$ be a fixed graph with $r$ edges, $s$ vertices and maximum degree
$\Delta$, and let $G=(V,E)$ be an $(n,d,\lambda)$-graph, where, say,
$d \leq 0.9n$.
Let $m <n$
satisfy $m \gg \lambda (\frac{n}{d})^{\Delta}$.
Then, for every
subset $V' \subset V$ of cardinality $m$, the number of (not necessarily
induced) copies of $H$ in $V'$ is
$$
\big(1+o(1)\big)\frac{m^s}{|Aut(H)|} \left(\frac{d}{n}\right)^r.
$$
\end{theo}
Note that this implies that a similar result holds for the number
of induced copies of
$H$. Indeed, if $ n \gg d$ and
$m \gg \lambda (\frac{n}{d})^{\Delta+1}$ then the number of copies
of each graph obtained from $H$ by adding to it at least one edge
is, by the above Theorem, negligible compared to the number of
copies of $H$, and hence almost all copies of $H$ in $V'$ are induced. If
$d=\Theta(n)$ then,
by inclusion-exclusion, the number of
induced copies of $H$ in $V'$ as above is also roughly the
"correct" number. A special case of the above theorem implies that if
$\lambda =O(\sqrt d)$ and $d \gg n^{2/3}$, then any $(n,d, \lambda)$-graph
contains many triangles. As shown in Example 9, Section 3, this is not true
when $d=(\frac{1}{4}+o(1)) n^{2/3}$, showing that the assertion of
the theorem is not far from being best possible.

\noindent
{\bf Proof of Theorem \ref{number-subgraphs}.}\, 
To prove the theorem, consider a random one-to-one mapping
of the set of vertices of $H$ into the set of vertices $V'$.
Denote by $A(H)$ the event that
every edge of $H$ is mapped on an edge of $G$. In such a
case we say that the mapping is an embedding of $H$. Note that it
suffices to prove that
\begin{equation}
\label{e1}
Pr(A(H))=(1+o(1))\left(\frac{d}{n}\right)^r.
\end{equation}

We prove (\ref{e1}) by induction on the number of edges $r$.  
The base case $(r=0)$ is trivial. Suppose that (\ref{e1}) holds for
all graphs with less than $r$ edges, and let
$uv$ be an edge of $H$.
Let $H_{uv}$ be the graph obtained from $H$ by removing the edge $uv$
(and keeping all vertices).
Let $H_u$ and $H_v$ be the induced subgraphs of $H$ on
the sets of vertices $V(H)\setminus \{ v\}$  and $V(H)\setminus
\{ u\}$, respectively, and let $H'$ be the induced subgraph of $H$ on
the set of
vertices $V(H)\setminus
\{ u,v\}$. Let $r'$ be the number of
edges of $H'$ and note that
$r-r' \leq 2(\Delta-1)+1=2\Delta -1. $
Clearly $Pr(A(H_{uv}))=Pr(A(H_{uv})|A(H'))\cdot
Pr(A(H'))$. Thus, by the induction hypothesis applied to
$H_{uv}$ and to $H'$:
\[
Pr(A(H_{uv})|A(H'))=(1+o(1)) \left(\frac{d}{n}\right)^{r-1-r'}.
\]
For an embedding $f'$ of $H'$, let $\nu(u,f')$ be the number
of extensions of $f'$ to an embedding of $H_u$ in $V'$; $\nu(v,f')$
denotes the same for $v$. Clearly, the number of extensions of
$f'$ to an embedding of $H_{uv}$ in $V'$ is at least
$\nu(u,f')\nu(v,f')-\min(\nu(u,f'),\nu(v,f'))$ and at most
$\nu(u,f')\nu(v,f')$. Thus we have
\[
\frac{\nu(u,f')\nu(v,f')-\min(\nu(u,f'),\nu(v,f'))}{(m-s+2)(m-s+1)}\leq
Pr\big(A(H_{uv})|f'\big)\leq
\frac{\nu(u,f')\nu(v,f')}{(m-s+2)(m-s+1)}.
\]
Taking expectation over all embeddings $f'$ the middle term
becomes $Pr(A(H_{uv})|A(H'))$, which is $(1+o(1))
(\frac{d}{n})^{r-1-r'}$. Note that by our choice of the
parameters and the well known fact that $\lambda =\Omega(\sqrt d)$,
the expectation of
the term $\min(\nu(u,f'),\nu(v,f'))~(~ \leq m)$ is negligible and we get
\[
E_{f'}\big(\nu(u,f')\nu(v,f')|\ A(H')\big)=(1+o(1)) m^2 \left(\frac{d}{n}\right)^{r-1-r'}.
\]
Now let $f$ be a random
one-to-one mapping of $V(H)$ into $V'$.
Let $f'$ be a fixed embedding of $H'$. Then
\[
Pr_f\big(A(H)|\ f|_{V(H)\setminus \{ u,v\}}=f'\big)=
\left(\frac{d}{n}\right)\frac{\nu(u,f')\nu(v,f')}{(m-s+2)(m-s+1)}+\delta,
\]
where $|\delta|\leq
\lambda\frac{\sqrt{\nu(u,f')\nu(v,f')}}{(m-s+2)(m-s+1)}$.
This follows from Theorem \ref{eigen}, where we take the possible
images of $u$ as the set $U$ and the possible images of $v$ as
the set $W$. Averaging over embeddings $f'$ we get
$Pr(A(H)|A(H'))$ on the left hand side. On the right hand side
we get $(1+o(1)) (\frac{d}{n})^{r-r'}$ from the first term plus the
expectation of the error term $\delta$. By Jensen's inequality,
the absolute value of this expectation is
bounded by
\[
\lambda\frac{\sqrt{E(\nu(u,f')\nu(v,f'))}}{(m-s+2)(m-s+1)}
=(1+o(1)) \frac{\lambda}{m} \left(\frac{d}{n}\right)^{(r-r'-1)/2}.
\]
Our assumptions on the parameters imply that this is negligible
with respect to the main term. Therefore
$ Pr(A(H))=Pr(A(H)|A(H')) \cdot Pr(A(H'))=(1+o(1)) \left(\frac{d}{n}\right)^r$,
completing the proof of
Theorem \ref{number-subgraphs}.
\hfill $\Box$

If we are only interested in the existence of one copy of $H$ then one can 
sometimes improve the conditions on $d$ and $\lambda$ in Theorem 
\ref{number-subgraphs}. For example if 
$H$ is a complete graph of order $r$ then the following result was proved in 
\cite{AK}.
\begin{prop}
\label{cliques}\cite{AK}
Let $G$ be an $(n,d,\lambda)$-graph. Then for every integer $r \geq 2$ every set of vertices 
of $G$ of size 
more than
$$\frac{(\lambda+1)n}{d}\left(1+\frac{n}{d}+\ldots+\Big(\frac{n}{d}\Big)^{r-2}\right)$$
contains a copy of a complete graph $K_r$.
\end{prop}

\noindent
In particular, when $d\geq\Omega(n^{2/3})$ and  $\lambda \leq O(\sqrt{d})$ then 
any $(n,d,\lambda)$-graph  contains a
triangle and as shows Example 9 in Section 3 this is tight. Unfortunately 
we do not know if 
this bound is also tight for $r \geq 4$. It would be interesting to construct
examples of
$(n,d,\lambda)$-graphs with $d=\Theta\big(n^{1-1/(2r-3)}\big)$ and $\lambda \leq 
O(\sqrt{d})$ which contain no copy of $K_r$. 

Finally we present one additional result about the existence of odd
cycles in pseudo-random graphs.

\begin{prop}
\label{cycles}
Let $k\geq 1$ be an integer and let $G$ be an $(n,d,\lambda)$-graph such that 
$d^{2k}/n\gg \lambda^{2k-1}$. Then $G$ contains a cycle of length $2k+1$.
\end{prop}

\noindent
{\bf Proof.}\, 
Suppose that $G$ contains no cycle of length $2k+1$.
For every two vertices $u,v$ of $G$ denote by $d(u,v)$ the length of a shortest
path from $u$ to $v$. For every $i \geq 1$ let $N_i(v)=\{u~|~d(u,v)=i\}$ be the set of 
all vertices in $G$ which are at distance exactly $i$ from $v$. In \cite{EFRS} Erd\H{o}s et 
al. proved that if $G$ contains no cycle of length 
$2k+1$ then for any $1 \leq i \leq k$ the induced graph $G[N_i(v)]$ contains 
an independent set of size $|N_i(v)|/(2k-1)$. This result together with 
Proposition \ref{ind-set} implies that for every vertex $v$  and for every $1 \leq i \leq 
k$, $|N_i(v)| \leq (2k-1)\lambda n/d$. Since $d^{2k}/n\gg \lambda^{2k-1}$ we have that
$\lambda=o(d)$. Therefore by Theorem \ref{eigen}
$$e\big(N_i(v)\big) \leq \frac{d}{2n}|N_i(v)|^2+\lambda |N_i(v)| \leq 
\frac{d}{n} \frac{(2k-1)\lambda n}{2d}|N_i(v)|+\lambda |N_i(v)|
<2k\lambda|N_i(v)|=o\big(d|N_i(v)|\big).$$

Next we prove by induction that for every $1\le i \leq k$, 
$\frac{|N_{i+1}(v)|}{|N_i(v)|}\geq (1-o(1))d^2/\lambda^2$. 
By the above discussion the number of edges spanned by 
$N_1(v)$ is $o(d^2)$ and therefore 
$e\big(N_1(v),N_2(v)\big)=d^2-o(d^2)=(1-o(1))d^2$.
On the other hand, by Theorem \ref{eigen} 
\begin{eqnarray*}
e\big(N_1(v),N_2(v)\big) &\leq& \frac{d}{n}|N_1(v)||N_2(v)|+\lambda\sqrt{|N_1(v)||N_2(v)|} 
\leq \frac{d}{n}\,d\,\frac{(2k-1)\lambda n}{d} +\lambda\sqrt{d|N_2(v)|}\\
&=&\lambda d\sqrt{\frac{|N_2(v)|}{d}}+O(\lambda d)=
\lambda d\sqrt{\frac{|N_2(v)|}{|N_1(v)|}}+o(d^2).
\end{eqnarray*}
Therefore $\frac{|N_2(v)|}{|N_1(v)|}\geq (1-o(1))d^2/\lambda^2$. Now
assume that
$\frac{|N_{i}(v)|}{|N_{i-1}(v)|}\geq (1-o(1))d^2/\lambda^2$. Since the number 
of edges spanned by $N_{i}(v)$ is $o\big(d|N_i(v)|\big)$ we obtain
\begin{eqnarray*}
e\big(N_i(v),N_{i+1}(v)\big)&=&d|N_i(v)|-2e\big(N_i(v)\big)-e\big(N_{i-1}(v),N_i(v)\big) \\
&\geq& d|N_i(v)|-o\big(d|N_i(v)|\big)-d|N_{i-1}(v)|\\
&\geq& 
(1-o(1))d|N_i(v)|-(1+o(1))d(\lambda^2/d^2)|N_i(v)|\\
&=&(1-o(1))d|N_i(v)|-o\big(d|N_i(v)|\big)=
(1-o(1))d|N_i(v)|.
\end{eqnarray*}
On the other hand, by Theorem \ref{eigen}
\begin{eqnarray*}
e\big(N_i(v),N_{i+1}(v)\big) &\leq& 
\frac{d}{n}|N_i(v)||N_{i+1}(v)|+\lambda\sqrt{|N_i(v)||N_{i+1}(v)|}\\
&\leq& \frac{d}{n} \,\frac{(2k-1)\lambda n}{d}\,|N_i(v)| 
+\lambda\sqrt{|N_i(v)||N_{i+1}(v)|}\\
&=&O(\lambda |N_i(v)|)+\lambda |N_i(v)|\sqrt{\frac{|N_{i+1}(v)|}{|N_i(v)|}}=
\lambda |N_i(v)|\sqrt{\frac{|N_{i+1}(v)|}{|N_i(v)|}}+o\big(d|N_i(v)|\big).
\end{eqnarray*}
Therefore $\frac{|N_{i+1}(v)|}{|N_i(v)|}\geq (1-o(1))d^2/\lambda^2$ and we proved the 
induction step.

Finally note that 
$$|N_k(v)| = d \prod_{i=1}^{k-1}\frac{|N_{i+1}(v)|}{|N_i(v)|}\geq
(1+o(1))d\left(\frac{d^2}{\lambda^2}\right)^{k-1}=
(1+o(1)) \frac{d^{2k-1}}{\lambda^{2k-2}}\gg (2k-1)\frac{\lambda n}{d}.$$
This contradiction completes the proof.
\hfill $\Box$

\vspace{0.15cm}
\noindent
This result implies that when $d\gg n^{\frac{2}{2k+1}}$ and $\lambda \leq O(\sqrt{d})$
then any $(n,d,\lambda)$-graph contains a cycle of length $2k+1$. As shown
by Example 10 of the 
previous section this result is tight. It is worth mentioning here that
it follows from the result of 
Bondy and Simonovits \cite{BS} that any $d$-regular graph with 
$d \gg n^{1/k}$ contains a cycle of length $2k$.
Here we do not need to make any assumption about the second eigenvalue $\lambda$.
This bound is known to be tight for $k=2,3,5$ (see Examples 6,7, Section 3).

\subsection{Extremal properties}
Tur\'an's theorem \cite{T} is one of the fundamental results in
Extremal Graph Theory. It  states that among $n$-vertex graphs not 
containing a clique of size $t$ the complete $(t-1)$-partite 
graph with (almost) equal parts has the maximum number of edges. 
For two graphs $G$ and $H$ we define the Tur\'an number $ex(G,H)$ of $H$
in $G$, as the largest integer $e$, such that there is an $H$-free subgraph of
$G$ with $e$ edges. Obviously $ex (G,H) \leq |E(G)|$, where $E(G)$ denotes the
edge set of $G$.  Tur\'an's theorem, in an asymptotic form, can be restated as 
$$ex(K_n, K_t) = \left(\frac{t-2}{t-1}+o(1)\right){n\choose 2},$$
that is the largest $K_t$-free subgraph of $K_n$
contains  approximately $\frac{t-2}{t-1}$-fraction of its edges.
Here we would like to describe an extension of this result to 
$(n,d,\lambda)$-graphs.

For an arbitrary  graph $G$ on $n$ vertices it is
easy to give a lower bound on $ex (G,K_t)$ following Tur\'an's
construction. One can partition the vertex set of $G$ into $t-1$ 
parts such that the degree of each vertex within its own part
is at most $\frac{1}{t-1}$-times its degree in $G$. Thus the subgraph
consisting of the edges of $G$ connecting two different parts has at least a
$\frac{t-2}{t-1}$-fraction of the edges of $G$ and
is clearly $K_t$-free. We say that a graph (or rather a family of graphs)
is {\em $t$-Tur\'an} if this
trivial lower bound is essentially an  upper bound as well. More precisely,
$G$ is $t$-Tur\'an if $ex(G,K_t) = \big(\frac{t-2}{t-1} +o(1)\big) |E(G)|$.

It has been shown that for any fixed $t$, there is a
number $m(t,n)$ such that almost all graphs on $n$ vertices with
$m \ge m(t,n) $ edges are $t$-Tur\'an 
(see \cite{SV}, \cite{KohRodS} for the most recent estimate
for $m(t,n)$). However, these 
results are about random graphs and do not provide a 
deterministic sufficient condition for a graph to be $t$-Tur\'an. 
It appears that such a condition can be obtained by a simple assumption
about the spectrum of the graph. This was proved by 
Sudakov, Szab\'o and Vu in \cite{SSV}. They obtained the following result. 

\begin{theo}
\label{t1} \cite{SSV}
Let $t\geq 3$ be an integer and let $G=(V,E)$ be an
$(n,d,\lambda)$-graph. If $\lambda=o(d^{t-1}/n^{t-2})$ then
$$ex(G,K_t)=\left(\frac{t-2}{t-1}+o(1)\right)|E(G)|.$$
\end{theo}
Note that this theorem generalizes Tur\'an's theorem, as
the second eigenvalue of the complete graph $K_n$ is 1. 

Let us briefly discuss the sharpness of Theorem \ref{t1}. 
For $t=3$, one can show that its condition involving $n,d$ and $\lambda$ is 
asymptotically tight. Indeed,
in this case the above theorem
states that if $d^2/n\gg \lambda$, then one needs to delete about
half of the edges of $G$ to destroy all the triangles. On
the other hand, by taking the example of Alon
(Section \ref{examples}, Example 9) whose parameters are:
$d=\Theta(n^{2/3})$, $\lambda=\Theta(n^{1/3})$, and blowing it up (which
means replacing each vertex by an independent set of size $k$ and
connecting two vertices in the new graph if and only if the
corresponding vertices of $G$ are connected by an edge) we get a graph
$G(k)$ with the following properties:

\begin{center}
$|V(G(k))|=n_k=nk$;\,\, $G(k)$ is $d_k=dk$-regular;\,\, $G(k)$ is
triangle-free;\,\,
$\lambda(G(k))=k\lambda$\, and \,$\lambda(G(k))=\Omega\big(d_k^2/n_k\big)$.
\end{center}

\noindent
The above bound for the second eigenvalue of $G(k)$ can be obtained by
using well known results on the eigenvalues of the tensor product of
two matrices, see \cite{KriSudSza02} for more details.  
This construction implies that for $t=3$ and any
sensible degree $d$ the condition in Theorem \ref{t1} is not far
from being best possible.

\subsection{Factors and fractional factors}
Let $H$ be a fixed graph on $n$ vertices. We say that a graph $G$ on
$n$ vertices has an {\em $H$-factor} if $G$ contains $n/h$ vertex
disjoint copies of $H$. Of course, a trivial necessary condition for the
existence of an $H$-factor in $G$ is that $h$ divides $n$. For example,
if $H$ is just an edge $H=K_2$, then an $H$-factor is a perfect matching
in $G$.

One of the most important classes of graph embedding problems is to find
sufficient conditions for the existence of an $H$-factor in a graph $G$,
usually assuming that $H$ is fixed while the order $n$ of $G$ grows. In
many cases such conditions are formulated in terms of the minimum
degree of $G$. For example, the classical result of Hajnal and
Szemer\'edi \cite{HajSze70} asserts that if the minimum degree
$\delta(G)$ satisfies $\delta(G)\ge (1-\frac{1}{r})n$, then $G$
contains $\lfloor n/r\rfloor$ vertex disjoint copies of $K_r$. The
statement of this theorem is easily seen to be tight.

It turns our that pseudo-randomness allows in many cases to
significantly weaken sufficient conditions for $H$-factors and to obtain
results which fail to hold for general graphs of the same edge density.

Consider first the case of a constant edge density $p$. In this case the
celebrated Blow-up Lemma of Koml\'os, S\'ark\"ozy and Szemer\'edi
\cite{KomSarSze97} can
be used to show the existence of $H$-factors. In order to formulate the
Blow-up Lemma we need to introduce the notion of a super-regular pair.
Given $\epsilon>0$ and $0<p<1$, a bipartite graph $G$ with bipartition
$(V_1,V_2)$, $|V_1|=|V_2|=n$, is called {\em super 
$(p,\epsilon)$-regular} if
\begin{enumerate}
\item For all vertices $v\in V(G)$,
$$
(p-\epsilon)n\le d(v)\le (p+\epsilon)n\ ;
$$ 
\item For every pair of sets $(U,W)$, $U\subset V_1$, $W\subset V_2$,
$|U|,|W|\ge \epsilon n$,
$$
\left|\frac{e(U,W)}{|U||W|}-\frac{|E(G)|}{n^2}\right|\le \epsilon\ .
$$
\end{enumerate}

\begin{theo}\label{KSS}\cite{KomSarSze97}
For every choice of integers $r$ and $\Delta$ and a real $0<p<1$ there
exist an $\epsilon>0$ and an integer $n_0(\epsilon)$ such that the
following is true. Consider an $r$-partite graph $G$ with all partition
sets $V_1,\ldots,V_r$ of order $n>n_0$ and all ${r\choose 2}$ bipartite
subgraphs $G[V_i,V_j]$ super $(p,\epsilon)$-regular. Then for every
$r$-partite graph $H$ with maximum degree $\Delta(H)\le \Delta$ and all
partition sets $X_1,\ldots,X_r$ of order $n$, there exists an embedding
$f$ of $H$ into $G$ with each set $X_i$ mapped onto $V_i$,
$i=1,\ldots,r$.
\end{theo}
(The above version of the Blow-up Lemma, due to R\"odl and Ruci\'nski
\cite{RodRuc99}, is somewhat different from and yet equivalent to
the original formulation of Koml\'os et al. We use it here as it is
somewhat closer in spirit to the notion of pseudo-randomness).

The Blow-up Lemma is a very powerful embedding tool. Combined with
another "big cannon", the Szemer\'edi Regularity Lemma, it can be used
to obtain approximate versions of many of the most famous embedding
conjectures. We suggest the reader to consult a survey of Koml\'os
\cite{Kom99} for more details and discussions.

It is easy to show that if $G$ is an $(n,d,\lambda)$-graph with
$d=\Theta(n)$ and $\lambda=o(n)$, and $h$ divides $n$, then a random
partition of $V(G)$ into $h$ equal parts $V_1,\ldots,V_h$ produces
almost surely ${h\choose 2}$ super $(d/n,\epsilon)$-regular pairs. Thus
the Blow-up Lemma can be applied to the obtained $h$-partite subgraph of
$G$ and we get:
\begin{coro}\label{Hdense}
Let $G$ be an $(n,d,\lambda)$-graph with $d=\Theta(n)$, $\lambda=o(n)$.
If $h$ divides $n$, then $G$ contains an $H$-factor, for every fixed
graph $H$ on $h$ vertices.
\end{coro}

The case of a vanishing edge density $p=o(1)$ is as usual
significantly more complicated. Here a sufficient condition for the
existence of an $H$-factor should depend heavily on the graph $H$, as
there may exist quite dense pseudo-random graphs without a single copy
of $H$, see, for example, the Alon graph (Example 9 of Section
\ref{examples}). When $H=K_2$, already a very weak pseudo-randomness
condition suffices to guarantee an $H$-factor, or a perfect matching, as
provided by Theorem \ref{edge-connectivity}. We thus consider the case
$H=K_3$, the task here is to guarantee a {\em triangle factor}, i.e. a
collection of $n/3$ vertex disjoint triangles. This problem has been
treated by Krivelevich, Sudakov and Szab\'o \cite{KriSudSza02}
who obtained the following result:
\begin{theo}\cite{KriSudSza02}\label{KrSS}
Let $G$ be an $(n,d,\lambda)$-graph. If $n$ is divisible by 3 and 
$$
\lambda=o\left(\frac{d^3}{n^2\log n}\right)\,,
$$
then $G$ has a triangle factor.
\end{theo}
For best pseudo-random graphs with $\lambda=\Theta(\sqrt{d})$ the
condition of the above theorem is fulfilled when $d\gg
n^{4/5}\log^{2/5}n$.

To prove Theorem \ref{KrSS} Krivelevich et al. first partition the
vertex set $V(G)$ into three parts $V_1,V_2,V_3$ of equal cardinality at
random. Then they choose a perfect matching $M$ between $V_1$ an $V_2$
at random and form an auxiliary bipartite graph $\Gamma$ whose parts
are $M$ and $V_3$, and whose edges are formed by connecting $e\in M$ and
$v\in V_3$ if both endpoints of $e$ are connected by edges to $v$ in
$G$. The existence of a perfect matching in $\Gamma$ is equivalent to
the existence of a triangle factor in $G$. The authors of
\cite{KriSudSza02} then proceed to show that if $M$ is chosen at random
then the Hall condition is satisfied for $\Gamma$ with positive
probability.

The result of Theorem \ref{KrSS} is probably not tight. In fact, the
following conjecture is stated in \cite{KriSudSza02}:
\begin{conj}\cite{KriSudSza02}\label{KrSScon}
There exists an absolute constant $c>0$ so that every $d$-regular graph
$G$ on $3n$ vertices, satisfying $\lambda(G)\le cd^2/n$, has a triangle
factor.
\end{conj} 
If true the above conjecture would be best possible, up to a constant
multiplicative factor. This is shown by taking the example of Alon
(Section \ref{examples}, Example 9) 
and blowing each of its vertices by an independent set of size $k$. 
As we already discussed in the previous section
(see also \cite{KriSudSza02}), this 
gives a triangle-free $d_k$-regular graph
$G(k)$ on $n_k$ vertices which satisfies
$\lambda(G(k))=\Omega\big(d_k^2/n_k\big)$.

Krivelevich, Sudakov and Szab\'o considered in \cite{KriSudSza02} also
the fractional version of the triangle factor problem.
 Given a graph $G=(V,E)$, denote by $T=T(G)$ the set of all triangles
of $G$. A function $f:T\rightarrow {\mathbb R}_+$ is called a 
{\em fractional triangle factor} if for every $v\in V(G)$ one has
$\sum_{v\in t}f(t)=1$. If $G$ contains a triangle factor $T_0$, then
assigning values $f(t)=1$ for all $t\in T_0$, and $f(t)=0$ for all other $t\in T$ produces a fractional
triangle factor. This simple argument shows that the existence of a
triangle factor in $G$ implies the existence of a fractional triangle
factor. The converse statement is easily seen to be invalid in general.

The fact that a fractional triangle factor $f$ can take non-integer
values, as opposed to the characteristic vector of a "usual" (i.e.
integer) triangle factor, enables to invoke the powerful machinery of
Linear Programming to prove a much better result than Theorem
\ref{KrSS}.
\begin{theo}\label{KrSSfr}\cite{KriSudSza02}
Let $G=(V,E)$ be a $(n,d,\lambda)$-graph If $\lambda \le
0.1d^2/n$ then $G$ has a fractional triangle factor.
\end{theo}
This statement is optimal up to a constant factor -- see the discussion
following Conjecture \ref{KrSScon}.

Already for the next case $H=K_4$ analogs of Theorem \ref{KrSS} and
\ref{KrSSfr} are not known. In fact, even an analog of Conjecture
\ref{KrSScon} is not available either, mainly due to the fact that we do not
know the weakest possible spectral condition guaranteeing a single copy
of $K_4$, or $K_r$ in general, for $r\ge 4$. 

Finally it would be interesting to show that
for every integer $\Delta$ there exist a real $M$ and an integer $n_0$
so that the following is true. If $n\ge n_0$ and $G$ is an
$(n,d,\lambda)$-graph for which 
$\lambda \le d\big(d/n\big)^M,$
then $G$ contains a copy of any graph $H$ on at most $n$ vertices with
maximum degree $\Delta(H)\le \Delta$. This can be considered as a
sparse analog of the Blow-up Lemma.

\subsection{Hamiltonicity}

A {\em Hamilton cycle} in a graph is a cycle passing through all the
vertices of this graph. A graph is called {\em Hamiltonian} if it has at
least one Hamilton cycle. For background information on Hamiltonian
cycles the reader can consult a survey of Chv\'atal \cite{Chv85}.

The notion of Hamilton cycles is one of the most central in modern Graph
Theory, and many efforts have been devoted to obtain sufficient
conditions for Hamiltonicity. The absolute majority of such known 
conditions
(for example, the famous theorem of Dirac asserting that a graph on $n$
vertices with minimal degree at least $n/2$ is Hamiltonian) deal with
graphs which are fairly dense. Apparently there are very few sufficient
conditions for the existence of a Hamilton cycle in sparse graphs.

As it turns out spectral properties of graphs can supply rather powerful
sufficient conditions for Hamiltonicity. Here is one such result, quite
general and yet very simple to prove, given our knowledge of properties
of pseudo-random graphs.

\begin{prop}\label{Ham1}
Let $G$ be an $(n,d,\lambda)$-graph. If
$$
d-36\frac{\lambda^2}{d} \ge \frac{\lambda n}{d+\lambda}\,,
$$
then $G$ is Hamiltonian.
\end{prop}

\noindent{\bf Proof.\ } According to Theorem \ref{connectivity} $G$ is
$(d-36\lambda^2/d)$-vertex-connected. Also, $\alpha(G)\le \lambda 
n/(d+\lambda)$, as 
stated 
in Proposition \ref{ind-set}. 
Finally, a theorem of Chv\'atal and Erd\H os \cite{ChvErd72}
asserts that if the vertex-connectivity of a graph $G$ is at least as 
large as its independence number, then $G$ is Hamiltonian.\hfill $\Box$

\medskip  

The Chv\'atal-Erd\H os Theorem has also been used by Thomason in
\cite{Tho87a}, who proved that a $(p,\alpha)$-jumbled graph $G$ with
minimal degree $\delta(G)=\Omega(\alpha/p)$ is Hamiltonian. His proof is
quite similar in spirit to that of the above proposition.

Assuming that $\lambda=o(d)$ and $d\rightarrow\infty$, the condition of
Proposition \ref{Ham1} reads then as: $\lambda\le (1-o(1))d^2/n$. For
best possible pseudo-random graphs, where $\lambda=\Theta(\sqrt{d})$,
this condition starts working when $d=\Omega(n^{2/3})$. 

One can however prove a much stronger asymptotical result, using more
sophisticated tools for assuring Hamiltonicity. The authors prove such a
result in \cite{KriSud02}:

\begin{theo}\label{Ham2}\cite{KriSud02}
Let $G$ be an $(n,d,\lambda)$-graph.  If $n$ is large enough and
$$
\lambda \leq \frac{(\log\log n)^2}{1000\log n(\log\log\log n)}\,d\,,
$$
then $G$ is Hamiltonian.
\end{theo}

The proof of Theorem \ref{Ham2} is quite involved technically. Its main
instrument is the famous rotation-extension technique of Posa
\cite{Pos76}, or rather a version of it developed by Koml\'os and
Szemer\'edi in \cite{KomSze83} to obtain the exact threshold for the
appearance of a Hamilton cycle in the random graph $G(n,p)$. We omit the
proof details here, referring the reader to \cite{KriSud02}.

For reasonably good pseudo-random graphs, in which
$\lambda\le d^{1-\epsilon}$ for some $\epsilon>0$, Theorem \ref{Ham2}
starts working already when the degree $d$ is only polylogarithmic in
$n$ -- quite a progress compared to the easy Proposition \ref{Ham1}! It
is possible though that an even stronger result is true as given by the 
following conjecture:

\begin{conj}\label{cHam}\cite{KriSud02}
There exists a positive constant $C$ such that for large enough $n$, any
$(n,d,\lambda)$-graph that satisfies $d/ \lambda>C$
contains a Hamilton cycle.
\end{conj}

This conjecture is closely related to another well known problem on
Hamiltonicity.  The {\em toughness} $t(G)$ of a graph $G$ is the largest
real $t$ so that for every positive integer $x \geq 2$
one should delete at least $tx$ vertices
from $G$ in order to get an induced subgraph of it with at least
$x$ connected components. $G$ is $t$-tough if $t(G) \geq t$.
This parameter was introduced by
Chv\'atal in \cite{Chv73}, where he observed that
Hamiltonian graphs are $1$-tough and conjectured that
$t$-tough graphs are Hamiltonian for large enough $t$.
Alon showed in \cite{Alo95} that if $G$ is an $(n,d,\lambda)$-graph,
then the toughness of $G$ satisfies $t(G)>\Omega(d/\lambda)$. Therefore
the conjecture of Chv\'atal implies the above conjecture.

Krivelevich and Sudakov used Theorem \ref{Ham2} in \cite{KriSud02} to
derive Hamiltonicity of sparse random Cayley graphs. Given a group $G$
of order $n$, choose a set $S$ of $s$ non-identity elements uniformly
at random and form a Cayley graph $\Gamma(G,S\cup S^{-1})$ (see 
Example 8 in Section 3 for the definition of a Cayley graph). The question 
is how large
should be the value of $t=t(n)$ so as to guarantee the almost sure
Hamiltonicity of the random Cayley graph no matter which group $G$ we
started with.

\begin{theo}\label{Ham3}\cite{KriSud02}
Let $G$ be a group of order $n$.
Then for every $c>0$ and large enough $n$
a Cayley graph $X(G,S\cup S^{-1})$, formed by
choosing a set $S$ of $c\log^5 n$ random generators in $G$,
is almost surely Hamiltonian.
\end{theo}

\noindent{\bf Sketch of proof.\ } 
Let $\lambda$ be the second largest by absolute value eigenvalue of
$X(G,S)$. Note that the Cayley graph $X(G,S)$ is $d$-regular for
$d \geq c\log^5 n$. Therefore to prove Hamiltonicity of
$X(G,S)$, by Theorem \ref{Ham2} it is enough to show that almost surely
$\lambda/d \leq O(\log n)$. This can be done by 
applying an approach of Alon and Roichman \cite{AloRoi94} for bounding
the second eigenvalue of a random Cayley graph.\hfill$\Box$

\medskip

We note that a well known conjecture claims that every connected Cayley
graph is Hamiltonian. If true the conjecture would easily imply that
as few as $O(\log n)$ random generators are enough to give almost sure
connectivity and thus Hamiltonicity.

\subsection{Random subgraphs of pseudo-random graphs}
There is a clear tendency in recent years to study random graphs
different from the classical by now model $G(n,p)$ of binomial random
graphs. One of the most natural models for random graphs, directly
generalizing $G(n,p)$, is defined as follows. Let $G=(V,E)$ be a graph
and let $0<p<1$. The {\em random subgraph} $G_p$ if formed by choosing
every edge of $G$ independently and with probability $p$. Thus, when $G$
is the complete graph $K_n$ we get back the probability space $G(n,p)$.
In many cases the obtained random graph $G_p$ has many interesting and
peculiar features, sometimes reminiscent of those of $G(n,p)$, and
sometimes inherited from those of the host graph $G$. 

In this subsection we report on various results obtained on random
subgraphs of pseudo-random graphs. While studying this subject, we
study in fact not a single probability space, but rather a family
of probability  spaces, having many common features, guaranteed by
those of pseudo-random graphs. Although several results have already
been achieved in this direction, overall it is much less developed than
the study of binomial random graphs $G(n,p)$, and one can certainly
expect many new results on this topic to appear in the future. 

We start with Hamiltonicity of random subgraphs of pseudo-random
graphs. As we learned in the previous section spectral condition are in
many cases sufficient to guarantee Hamiltonicity. Suppose then that a
host graph $G$ is a Hamiltonian $(n,d,\lambda)$-graph. How small can 
the edge probability  $p=p(n)$ be chosen so as to guarantee almost sure
Hamiltonicity of the random subgraph $G_p$? This question has been
studied by Frieze and the first author in \cite{FriKri02}. They obtained
the following result.

\begin{theo}\label{FK}\cite{FriKri02}
Let $G$ be an $(n,d,\lambda)$-graph.
Assume that
$\lambda=o\left(\frac{d^{5/2}}{n^{3/2}(\log n)^{3/2}}\right)$. Form a
random subgraph $G_p$ of $G$ by choosing each edge of $G$ independently
with probability $p$. Then for any function $\om(n)$ tending to infinity
arbitrarily slowly:
\begin{enumerate}
\item if $p(n)=\frac{1}{d}(\log n+\log\log n-\om(n))$, then 
$G_p$ is almost surely not Hamiltonian;
\item if $p(n)=\frac{1}{d}(\log n+\log\log n+\om(n))$, then 
 $G_p$ is almost surely Hamiltonian.
\end{enumerate}
\end{theo} 

Just as in the case of $G(n,p)$ (see, e.g. \cite{Bol01}) it is quite
easy to predict the critical probability for the appearance of a
Hamilton cycle in $G_p$. An obvious obstacle for its existence is a
vertex of degree at most one. If such a vertex almost surely exists in
$G_p$, then $G_p$ is almost surely non-Hamiltonian. It is a
straightforward exercise to show that the smaller probability in the
statement of Theorem \ref{FK} gives the almost sure existence of such
a vertex. The larger probability can be shown to be sufficient to
eliminate almost surely all vertices of degree at most one in $G_p$.
Proving that this is sufficient for almost sure Hamiltonicity is much
harder. Again as in the case of $G(n,p)$ the rotation-extension
technique of Posa \cite{Pos76} comes to our rescue. We omit technical
details of the proof of Theorem \ref{FK}, referring the reader to
\cite{FriKri02}.

One of the most important events in the study of random graphs was the
discovery of the sudden appearance of the giant component by Erd\H os
and R\'enyi \cite{ErdRen60}. They proved that all connected components
of $G(n,c/n)$ with $0<c<1$ are almost surely trees or unicyclic and
have size $O(\log n)$. On the other hand, if $c>1$, then $G(n,c/n)$
contains almost surely a unique component of size linear in $n$ (the so
called {\em giant component}), while all other components are at most
logarithmic in size. Thus, the random graph $G(n,p)$ experiences the so
called {\em phase transition} at $p=1/n$. 

Very recently Frieze, Krivelevich and Martin showed \cite{FriKriMar02}
that a very similar behavior holds for random subgraphs of many
pseudo-random graphs. To formulate their result, for $\alpha>1$ we
define  $\bar{\alpha}<1$ to be the unique solution (other than 
$\alpha$) of the equation $xe^{-x}=\alpha e^{-\alpha}$.
 
\begin{theo}\label{FKM}\cite{FriKriMar02}
Let $G$ be an $(n,d,\lambda)$-graph. Assume that 
$
\lambda=o(d).
$
Consider the random subgraph $G_{\alpha/d}$, formed by choosing each
edge of $G$ independently and with probability $p=\alpha/d$. Then:
\begin{itemize}
\item[(a)] If $\alpha<1$ then almost surely the maximum component size 
is $O(\log n)$.
\item[(b)] If $\alpha>1$ then almost surely there is a unique giant
component of asymptotic size 
$\left(1-\frac{\bar{\alpha}}{\alpha}\right)n$
and the remaining components are of size $O(\log n)$.
\end{itemize}
\end{theo}

Let us outline briefly the proof of Theorem \ref{FKM}. First, bound
(\ref{eig}) and known estimates on the number of $k$-vertex trees in
$d$-regular graphs are  used to get estimates on the expectation of the
number of connected components of size $k$ in $G_{p}$, for various
values of $k$. Using these estimates it is proved then that almost
surely $G_p$ has no connected components of size between
$(1/\alpha\gamma)\log n$ and $\gamma n$ for a properly chosen
$\gamma=\gamma(\alpha)$. Define $f(\alpha)$ to be 1 for all $\alpha\le
1$, and to be $\bar{\alpha}/\alpha$ for $\alpha>1$. One can show then
that almost surely in $G_{\alpha/d}$
the number of vertices in components of size between 1 and
$d^{1/3}$ is equal to $nf(\alpha)$ up to the error term which is
$O(n^{5/6}\log n)$. This is done by first calculating the expectation of
the last quantity, which is asymptotically equal to $nf(\alpha)$, and
then by applying the Azuma-Hoeffding martingale inequality.    

Given the above, the proof of Theorem \ref{FKM} is straightforward. For
the case $\alpha<1$ we have $nf(\alpha)=n$ and therefore all but at
most $n^{5/6}\log n$ vertices lie in components of size at most
$(1/\alpha\gamma)\log n$. The remaining vertices should be in components
of size at least $\gamma n$, but there is no room for such components.
If $\alpha>1$, then $(\bar{\alpha}/\alpha)n+O(n^{5/6}\log n)$ vertices
belong to components of size at most $(1/\alpha\gamma)\log n$, and all
remaining vertices are in components of size at least $\gamma n$. These
components are easily shown to merge quickly into one giant component
of a linear size. The detail can be found in \cite{FriKriMar02} 
(see also \cite{ABS} for some related results).

One of the recent most popular subjects in the study of random graphs is
proving sharpness of thresholds for various combinatorial properties. 
This direction of research was spurred by a powerful theorem of
Friedgut-Bourgain \cite{Fri99}, providing a sufficient condition for the
sharpness of a threshold. The authors together with Vu apply this
theorem in \cite{KriSudVu02}
to show sharpness of graph connectivity, sometimes also called
{\em network reliability}, in random subgraphs of a wide class of graphs.
Here are the relevant definitions. For a connected graph $G$ and edge
probability $p$ denote by $f(p)=f(G,p)$ the probability that a random
subgraph $G_p$ is connected. The function $f(p)$ can be easily shown to
be strictly monotone.  For a fixed
positive constant $x \leq 1$ and a graph $G$, let
$p_{x} $ denote the (unique) value of $p$ where $f(G, p_{x})=
x$. We say that a family $(G_i)_{i=1}^{\infty}$ of graphs
satisfies the {\em sharp threshold} property if for any fixed positive
$\epsilon \le 1/2$
$$
\lim_{i \rightarrow \infty} \frac{ p_{\epsilon} (G_i)} { p_{1-\epsilon}
(G_i)} \rightarrow 1. 
$$ 
Thus the threshold for connectivity is sharp if the width of the
transition interval is negligible compared to the critical probability.
Krivelevich, Sudakov and Vu proved in \cite{KriSudVu02} the following theorem.

\begin{theo}\cite{KriSudVu02}\label{KSV}
Let $(G_i)_{i=1}^{\infty}$ be a family of distinct
graphs, where $G_i$ has $n_i$ vertices, maximum degree $d_i$ and
it is $k_i$-edge-connected. If
$$
\lim_{i \rightarrow \infty} \frac{k_i \ln n_i}{d_i}=\infty,
$$
then the family $(G_i)_{i=1}^{\infty}$ has a sharp connectivity
threshold.
\end{theo}
The above theorem extends a celebrated result of Margulis \cite{Mar74}
on network reliability (Margulis' result applies to the case where the
critical probability is a constant). 

Since $(n,d,\lambda)$ graphs are $d(1-o(1))$-connected as long as
$\lambda=o(d)$ by Theorem \ref{connectivity}, we immediately get the
following result on the sharpness of the connectivity threshold for
pseudo-random graphs.

\begin{coro}\label{KSVpr}
Let $G$ be an $(n,d,\lambda)$-graph. If $\lambda=o(d)$, then the
threshold for connectivity in the random subgraph $G_p$ is sharp.
\end{coro} 

Thus already weak connectivity is sufficient to guarantee sharpness of
the threshold. This result has potential practical applications as
discussed in \cite{KriSudVu02}.

Finally we consider a different probability space created from a graph
$G=(V,E)$. This space is obtained by putting random weights on the edges
of $G$ independently. One can then ask about the behavior of optimal
solutions for various combinatorial optimization problems. 

Beveridge, Frieze and McDiarmid treated in \cite{BevFriMcD98} the
problem of estimating the weight of a random minimum length spanning
tree in regular graphs. For each edge $e$ of a connected
graph $G=(V,E)$ define the length $X_e$ of $e$ to be a random variable 
uniformly distributed
in the interval $(0,1)$, where all $X_e$ are independent. Let
$mst(G,\bf{X})$ denote the minimum length of a spanning tree in such a
graph, and let $mst(G)$ be the expected value of $mst(G,{\bf X})$.  Of
course, the value of $mst(G)$ depends on the connectivity structure of
the graph $G$. Beveridge et al. were able to prove however that if the
graph $G$ is assumed to be almost regular and has a modest edge
expansion, then $mst(G)$ can be calculated asymptotically:

\begin{theo}\cite{BevFriMcD98}\label{BFM}
Let $\alpha=\alpha(d)=O(d^{-1/3})$ and let $\rho(d)$ and $\omega(d)$
tend to infinity with $d$. Suppose that the graph $G=(V,E)$ satisfies 
$$
d\le d(v)\le (1+\alpha)d\quad\mbox{for all $v\in V(G)$}\,,
$$  
and
$$
\frac{e(S,V\setminus S)}{|S|}\ge \omega d^{2/3}\log d
\quad\mbox{for all $S\subset V$ with $d/2<|S|\le min\{\rho d,|V|/2\}$}\ .
$$
Then 
$$
mst(G)=(1+o(1))\frac{|V|}{d}\zeta(3)\ ,
$$
where the $o(1)$ term tends to 0 as $d\rightarrow\infty$, and
$\zeta(3)=\sum_{i=1}^{\infty}i^{-3}=1.202...$.
\end{theo}
The above theorem extends a celebrated result of Frieze \cite{Fri85},
who proved it in the case of the complete graph $G=K_n$. 

Pseudo-random graphs supply easily the degree of edge expansion required
by Theorem \ref{BFM}. We thus get:
\begin{coro}\label{BFMpr}
Let $G$ be an $(n,d,\lambda)$-graph. If $\lambda=o(d)$ then 
$$
mst(G)=(1+o(1))\frac{n}{d}\zeta(3)\ .
$$
\end{coro}

Beveridge, Frieze and McDiarmid also proved that the random variable
$mst(G,{\bf X})$ is sharply concentrated around its mean given by
Theorem \ref{BFM}.

Comparing between the very well developed research of binomial random 
graphs
$G(n,p)$ and few currently available results on random subgraphs of
pseudo-random graphs, we can say that many interesting problems
in the latter subject are yet to be addressed, such as the
asymptotic behavior of the independence number and the chromatic number,
connectivity, existence of matchings and factors, spectral properties,
to mention just a few.

\subsection{Enumerative aspects}
Pseudo-random graphs on $n$ vertices with edge density $p$ are quite
similar in many aspects to the random graph $G(n,p)$. One can thus
expect that counting statistics in pseudo-random graphs will be close to
those in truly random graphs of the same density. As the random graph
$G(n,p)$ is a product probability space in which each edge behaves
independently, computing the expected number of most subgraphs in
$G(n,p)$ is straightforward. Here are just a few examples:
\begin{itemize}
\item The expected number of perfect matchings in $G(n,p)$ is 
$\frac{n!}{(n/2)!2^{n/2}}p^{n/2}$ (assuming of course that $n$ is even);
\item The expected number of spanning trees in $G(n,p)$ is
$n^{n-2}p^{n-1}$;
\item The expected number of Hamilton cycles in $G(n,p)$ is 
$\frac{(n-1)!}{2}p^n$.
\end{itemize}
In certain cases it is possible to prove that the actual number of
subgraphs in a pseudo-random graph on $n$ vertices with edge density
$p=p(n)$ is close to the corresponding expected value in the binomial
random graph $G(n,p)$.

Frieze in \cite{Fri00} gave estimates on the number of perfect
matchings and Hamilton cycles in what he calls super
$\epsilon$-regular graphs. Let $G=(V,E)$ be a graph on $n$ vertices with
${n\choose 2}p$ edges, where $0<p<1$ is a constant. 
Then $G$ is called {\em super $(p,\epsilon)$-regular}, for a constant
$\epsilon>0$, if 
\begin{enumerate}
\item For all vertices $v\in V(G)$,
$$
(p-\epsilon)n \le d(v)\le (p+\epsilon)n\,;
$$
\item 
For all $U,W\subset V$, $U\cap W=\emptyset$, $|U|, |W|\ge \epsilon n$, 
$$
\left|\frac{e(U,W)}{|U||W|}-p\right|\le \epsilon\ .
$$
\end{enumerate}
Thus, a super $(p,\epsilon)$-regular graph $G$ can be considered a
non-bipartite analog of the notion of a super-regular pair defined
above.
In our terminology, $G$ is a weakly pseudo-random graph of constant
density $p$, in which {\em all} degrees are asymptotically equal to
$pn$. Assume that $n=2\nu$ is even. Let $m(G)$ denote the number of
perfect matchings in $G$ and  let $h(G)$ denote the number of Hamilton
cycles in $G$,
and let $t(G)$ denote the number of spanning trees in $G$.

\begin{theo}\label{F}\cite{Fri00}
If $\epsilon$ is sufficiently small and $n$ is sufficiently large then
\begin{description}
\item[{\bf (a)}]
$$ 
(p-2\epsilon)^{\nu}\frac{n!}{\nu!2^{\nu}}\le m(G) \le 
(p+2\epsilon)^{\nu}\frac{n!}{\nu!2^{\nu}}\ ;
$$
\item [{\bf (b)}]
$$
(p-2\epsilon)^nn!\le h(G)\le (p+2\epsilon)^nn!\ ;
$$
\end{description}
\end{theo}

Theorem \ref{F} thus implies that the numbers of perfect matchings and
of Hamilton cycles in super $\epsilon$-regular
graphs are quite close asymptotically to the expected values of the
corresponding quantities in the random graph $G(n,p)$. Part (b) of
Theorem \ref{F} improves significantly Corollary 2.9 of Thomason
\cite{Tho87a} which estimates from below the number of Hamilton cycles
in jumbled graphs.

Here is a very brief sketch of the proof of Theorem \ref{F}. To estimate
the number of perfect matchings in $G$, Frieze takes a random partition
of the vertices of $G$ into two equal parts $A$ and $B$ and estimates
the number of perfect matchings in the bipartite subgraph of $G$
between $A$ and $B$. This bipartite graph is almost surely super
$2\epsilon$-regular, which allows to apply bounds previously
obtained by Alon, R\"odl and Ruci\'nski \cite{AloRodRuc98} for such
graphs. 

Since each Hamilton cycle is a union of two perfect matchings, it
follows immediately that $h(G)\le m^2(G)/2$, establishing the desired
upper bound on $h(G)$. In order to prove a lower bound, let $f_k$ be the
number of 2-factors in $G$ containing exactly $k$ cycles, so that
$f_1=h(G)$. Let also $A$ be the number of ordered pairs of edge disjoint
perfect matchings in $G$. Then 
\begin{equation}\label{A}
A= \sum_{i=1}^{\lfloor n/3\rfloor}2^kf_k\ .
\end{equation}
For a perfect matching $M$ in $G$ let $a_M$ be the number of perfect
matchings of $G$ disjoint from $M$. Since deleting $M$ disturbs
$\epsilon$-regularity of $G$ only marginally, one can use part (a) of
the theorem to get $a_M\ge (p-2\epsilon)^{\nu}\frac{n!}{\nu!2^{\nu}}$.
Thus 
\begin{equation}\label{A1}
A=\sum_{M\in G}a_M\ge \left((p-2\epsilon)^{\nu}\frac{n!}{\nu!2^{\nu}}
\right)^2 \ge (p-2\epsilon)^nn!\cdot \frac{1}{3n^{1/2}}\ .
\end{equation}
Next Frieze shows that the ratio $f_{k+1}/f_k$ can be bounded by a
polynomial in $n$ for all $1\le k\le k_1=O(p^{-2})$, 
$f_k\le 5^{-(k-k_0)/2}\max\{f_{k_0+1},f_{k_0}\}$ for all $k\ge k_0+2,
k_0=\Theta(p^{-3}\log n)$ and that the ratio
$(f_{k_1+1}+\ldots+f_{\lfloor n/3\rfloor})/f_{k_1}$ is also bounded by
a polynomial in $n$.
Then from (\ref{A}), $A\le O_p(1)\sum_{k=1}^{k_0+1}f_k$ and thus $A\le
n^{O(1)}f_1$. Plugging (\ref{A1}) we get the desired lower bound. 

One can also show (see \cite {Noga}) that the number of spanning trees $t(G)$ in super
$(p,\epsilon)$-regular graphs satisfies:
$$
(p-2\epsilon)^{n-1}n^{n-2}\le t(G)\le (p+2\epsilon)^{n-1}n^{n-2}\ ,
$$
for small enough $\epsilon>0$ and large enough $n$.
In order to estimate from below the number of spanning trees in $G$,
consider a random mapping $f:V(G)\rightarrow V(G)$, defined by choosing
for each $v\in V$ its neighbor $f(v)$ at random. Each such $f$ defines
a digraph $D_f=(V,A_f)$, $A_f=\{(v,f(v)): v\in V\}$. Each component of
$D_f$ consists of cycle $C$ with a rooted forest whose roots are all in
$C$. Suppose that $D_f$ has $k_f$ components. Then a spanning tree
of $G$ can be obtained by deleting the lexicographically first edge of
each cycle in $D_f$, and then extending the $k_f$ components to a
spanning tree. Showing that $D_f$ has typically $O(\sqrt{n})$ 
components implies that most of the mappings $f$ create a digraph close
to a spanning tree of $G$, and therefore:
$$
t(G)\ge n^{-O(\sqrt{n})}|f:V\rightarrow V|\ge
n^{-O(\sqrt{n})}(p-\epsilon)n^n\ .
$$
For the upper bound on $t(G)$ let 
$\Omega^*=\{f:V\rightarrow V: (v,f(v))\in E(G)$ for $v\not =1$ and 
$f(1)=1 \}$. Then 
$$
t(G)\le |\Omega^*| \le \big((p+\epsilon)n\big)^{n-1}\le 
(p+2\epsilon)^{n-1}n^{n-2}\ .
$$
To see this consider the following injection from the spanning trees of
$G$ into $\Omega^*$: orient each edge of a tree $T$ towards vertex 1 and
set $f(1)=1$. Note that this proof does not use the fact that the graph is 
pseudo-random.
Surprisingly it shows that all nearly regular connected graphs with the 
same density have 
approximately the same number of spanning trees.

For sparse pseudo-random graphs one can use Theorem \ref{FK} to estimate
the number of Hamilton cycles. Let $G$ be an $(n,d,\lambda)$-graph
satisfying the conditions of Theorem \ref{FK}. Consider the random
subgraph $G_p$ of $G$, where $p=(\log n+2\log\log n)/d$. Let $X$ be the
random variable counting the number of Hamilton cycles in $G_p$.
According to Theorem \ref{FK}, $G_p$ has almost surely a Hamilton cycle,
and therefore $E[X]\ge 1-o(1)$. On the other hand, the probability that
a given Hamilton cycle of $G$ appears in $G_p$ is exactly $p^n$.
Therefore the linearity of expectation implies $E[X]=h(G)p^n$. Combining
the above two estimates we derive:
$$
h(G)\ge \frac{1-o(1)}{p^n}=\left(\frac{d}{(1+o(1))\log n}\right)^n
\ .
$$
We thus get the following corollary:

\begin{coro}\cite{FriKri02}\label{FKc}
 Let $G$ be an $(n,d,\lambda)$-graph with
$\lambda=o(d^{5/2}/(n^{3/2}(\log n)^{3/2}))$.
Then $G$ contains at least $\left(\frac{d}{(1+o(1))\log n}\right)^n$
Hamilton cycles.
\end{coro}

Note that the number of Hamilton cycles in any $d$-regular graph on $n$
vertices obviously does not exceed $d^n$. Thus for graphs satisfying the
conditions of Theorem \ref{FK} the above corollary provides an
asymptotically  tight estimate on the exponent of the number of Hamilton
cycles. 

\section{Conclusion}
Although we have made an effort to provide a systematic coverage of the
current research in pseudo-random graphs, there are certainly quite a 
few subjects that were left outside this survey, due to the limitations
of space and time (and of the authors' energy). Probably the most notable
omission is a discussion of diverse applications of pseudo-random
graphs to questions from other fields, mostly Extremal Graph Theory,
where pseudo-random graphs provide the best known bounds for an amazing
array of problems. We hope to cover this direction in one of our future
papers. Still, we would like to believe that this survey can be helpful
in mastering various results and techniques pertaining to this field. 
Undoubtedly many more of them are bound to appear in the future and will
make this fascinating subject even more deep, diverse and appealing.

\medskip

\noindent{\bf Acknowledgment.\ } The authors would like to thank Noga
Alon for many illuminating discussions and for kindly granting us his
permission to present his Theorem \ref{number-subgraphs} here. The 
proofs of Theorems \ref{connectivity}, \ref{edge-connectivity} were 
obtained in discussions with him.

\end{document}